\documentclass[11pt]{article}
\usepackage[american]{babel}
\usepackage{microtype}
\usepackage{graphicx}
\usepackage{float}
\usepackage{epstopdf}
\usepackage{pdfsync}
\usepackage{amsmath,amssymb,amsthm,amsfonts,mathrsfs}
\usepackage{amssymb}
\usepackage{verbatim}
\usepackage{setspace}
\usepackage{hyperref}
\usepackage[a4paper]{geometry}
\newtheorem{thm}{Theorem}[section]
\newtheorem{cor}[thm]{Corollary}
\newtheorem{lem}[thm]{Lemma}
\newtheorem{prop}[thm]{Proposition}

\theoremstyle{definition}

\theoremstyle{remark}

\numberwithin{equation}{section}

\DeclareMathSymbol{\C}{\mathalpha}{AMSb}{"43}

\newcommand{\lam}{\lambda}

\newcommand{\R}{{\mathbb{R}}}

\newcommand{\inte}{\int_{\mathbb{R}^d}}

\newcommand{\bsub}{\begin{subequations}}
\newcommand{\esub}{\end{subequations}$\!$}

\begin{document}
\title{Limit Behavior of Mass Critical Hartree Minimization Problems with Steep Potential Wells}
\author{
Yujin Guo\thanks{Wuhan Institute of Physics and Mathematics,
Chinese Academy of Sciences, P.O. Box 71010, Wuhan 430071,
    P. R. China.  Email: \texttt{yjguo@wipm.ac.cn}. Y. J.  Guo is partially supported by NSFC under Grant No. 11671394 and MOST under Grant No. 2017YFA0304500.
    },
    \, Yong Luo\thanks{University of Chinese Academy of Sciences, Beijing 100190, P. R. China;  Wuhan Institute of Physics and Mathematics,
    Chinese Academy of Sciences, P.O. Box 71010, Wuhan 430071,
    P. R. China.  Email: \texttt{luoyong.wipm@outlook.com}.  },
    \, and\, Zhi-Qiang Wang\thanks{Center for Applied Mathematics, Tianjin University, 300072 Tianjin, P. R. China;
Department of Mathematics and Statistics
Utah State University, Logan, UT 84322, USA.  Email: \texttt{zhi-qiang.wang@usu.edu}. Z.-Q. Wang is partially supported by NSFC under Grant No. 11771324.}
}

\date{\today}
\smallbreak \maketitle

\begin{abstract}
We consider  minimizers of the following mass critical Hartree minimization problem:
 \[ e_\lambda(N):=\underset{\{u\in H^1(\R^d),\,\|u\|^2_2=N\}}{\inf} E_\lambda(u),\,\ d\ge 3,
\]
where the Hartree energy functional $E_\lambda(u)$ is defined by
\[
 E_\lambda(u):=\int_{\R ^d}|\nabla u(x)|^2dx+\lambda \int_{\R ^d}g(x)u^2(x)dx-\frac{1}{2} \int_{\R ^d}\int_{\R ^d} \frac{u^2(x)u^2(y)}{|x-y|^2}dxdy,\,\ \lam >0,\]
and the steep potential $g(x)$ satisfies $0=g(0)=\inf _{\R^d}g(x)\le g(x)\le 1$ and $1-g(x)\in L^{\frac{d}{2}}(\R^d)$.
We prove that there exists a constant $N^*>0$, independent of $\lam g(x)$, such that if $N\ge N^*$, then $e_\lambda(N)$ does not admit minimizers for any $\lam >0$;  if  $0<N<N^*$, then there exists a constant $\lam ^*(N)>0$ such that $e_\lambda(N)$ admits minimizers for any $\lam >\lam ^*(N)$, and $e_\lambda(N)$ does not admit minimizers for $0<\lam <\lam ^*(N)$.   For any given $0<N<N^*$,  the limit behavior of positive minimizers for $e_\lambda(N)$ is also studied as $\lambda\to\infty$, where the mass concentrates at the bottom of $g(x)$.
\end{abstract}

\vskip 0.2truein

\section{Introduction}
In the 1980s, Lions et al. analyzed in \cite{CL,Lions1} the following Hartree minimization problem:
 \begin{equation}\label{A:Lions}
 \begin{split} \inf\Big\{&\int_{\R ^d}\big(|\nabla u(x)|^2+V(x)u^2(x)\big)dx-\frac{1}{2} \int_{\R ^d}\int_{\R ^d} \frac{u^2(x)u^2(y)}{|x-y| }dydx,\\
 &\,\ u\in H^1(\R^d),\, \ \int_{\R ^d}u^2dx=N>0\Big\},
\end{split}\end{equation}
where $d \ge 3 $ and $ V(x) \ge  0$ is an external potential. The problem (\ref{A:Lions}) arises in Quantum Mechanics, see \cite{EHL2,Lions1} for some discussions on the relevance of this problem to Physics. Especially, the existence of positive minimizers for (\ref{A:Lions}) was discussed in \cite{CL,Lions1} by developing and applying the celebrated concentration-compactness principle. After these pioneering works, the problem (\ref{A:Lions})
was studied widely over the past few decades. We remark that the problem (\ref{A:Lions}) is essentially a mass subcritical
problem, and some analysis of (\ref{A:Lions}) can be extended naturally to the general case where the term $\frac{1}{|x-y|}$ in (\ref{A:Lions}) is replaced by $\frac{1}{|x-y|^s}$ for all $0<s<2$.

On the other hand, the mass critical minimization problem, arising in Bose-Einstein condensates (BEC) in $\R^2$, was analyzed recently in   \cite{GDL,GS,GWZZ,GZZ} and references therein, where the authors focused on the class where $V(x)$  satisfies
 \begin{equation}\label{A:V}
 0\le V(x)\in L^\infty_{\rm loc}(\R^2)\ \ \mbox{and}\,\ \lim_{|x|\to\infty} V(x) = \infty.
\end{equation}
Under the assumption (\ref{A:V}), the mass critical case of (\ref{A:Lions}), where the term $\frac{1}{|x-y|}$  is replaced by $\frac{1}{|x-y|^2}$, was studied in \cite{DLS} and somewhere else.
Stimulated by above facts, in this paper we are interested in the mass critical case of (\ref{A:Lions}) with a steep potential well $V(x)\in L^\infty(\R^d)$. We remark that the steep potential wells were considered in \cite{BA,BW} in the setting of nonlinear Schr$\ddot{o}$dinger equations without constraints. And here we consider normalized solutions subject to $L^2$ norm constraints which are more in line of concern with stability issues of standing waves.

More precisely,  in this paper we investigate the following mass critical minimization problem:
 \begin{equation}\label{def:ea}
 e_\lambda(N):=\underset{u\in H^1(\R^d),\,\|u\|^2_2=N}{\inf} E_\lambda(u),\,\,d\geq3,\,
\end{equation}
where the Hartree energy functional $E_\lambda(u)$ satisfies
 \begin{equation}\label{def:Ea}
 E_\lambda(u):=\int_{\R ^d}|\nabla u(x)|^2dx+\lambda \int_{\R ^d}g(x)u^2(x)dx-\frac{1}{2} \int_{\R ^d}\int_{\R ^d} \frac{u^2(x)u^2(y)}{|x-y|^2}dxdy,\, \lambda>0.
 \end{equation}
Here we consider the steep potential $g(x)\in C^{\alpha}_{loc}(\R^d)$ with $\alpha \in (0,1)$  satisfying
\begin{itemize}
\item [$\rm(M_1).$] $0=g(0)=\underset{x\in \R^d}{\inf}g(x)\leq g(x)\le 1,\, 1-g(x)\in L^{\frac{d}{2}}(\R^d)$, where $0\in \mathbb{R}^{d}$ is the unique global minimum point of $g(x)$ in $\R^d$.
\end{itemize}
The main purposes of this paper are to classify  the existence  and nonexistence of minimizers for $e_\lambda(N)$, based on which we shall investigate the limit behavior of minimizers  as $\lam\to\infty$. Even though there are existing papers as mentioned above of studying mass critical minimization problems,
as far as we know, this paper might be the first work of studying mass critical minimization problems under the steep potential $g(x)$ satisfying $(M_1)$.
More importantly, our analysis shows that there appear new and interesting phenomena on the existence and nonexistence of minimizers for $e_\lambda(N)$, and the limit behavior of minimizers as $\lam\to\infty$ presents new and challenging difficulties, for which one needs to investigate new analytic approaches.

Related to the minimization problem $ e_\lambda(N)$, to state our main results we now introduce the following nonlocal Hartree equation
\begin{equation}\label{eq:He}
-\triangle u+u-\Big(\int_{\R^d}\frac{u^2(y)}{|x-y|^2}dy\Big)u=0\ \  \hbox{in}\,\   \R^d,\ \,0<u\in H^1(\R^d)\,\  \hbox{and}\,\ d\ge 3.
\end{equation}
We define the energy  functional of (\ref{eq:He}) by
\begin{equation}\label{eq:12}
I(u):=\frac{1}{2}\int_{\R ^d}|\nabla u(x)|^2+u^2(x)dx-\frac{1}{4} \int_{\R ^d}\int_{\R ^d} \frac{u^2(x)u^2(y)}{|x-y|^2}dxdy,\  u\in H^1(\R^d).
\end{equation}
Consider
\begin{equation}\label{eq:13}
S:=\Big\{u(x)\in H^1(\R^d):\, u(x)\ \hbox{is a positive solution of (\ref{eq:He})}\Big\},
\end{equation}
and
\begin{equation}\label{eq:17}
  G:=\Big\{u(x)\in S: \, I(u)\leq I(v)\hbox{ for all}\  v\in S\Big\}.
\end{equation}
We then say that $0<u\in H^1(\R^d)$ is a ground state of (\ref{eq:He}) if $u\in G$. Note from  \cite{LMLZ,MV} that  (\ref{eq:He}) admits ground states $Q>0$, which must be radially symmetric and admit the following exponential decay
\begin{equation}\label{eq:ed}
 Q(|x|),\, |\nabla Q(|x|)|=O(e^{-\mu |x|})\ \,\, \hbox{as}\,\,\, |x|\rightarrow\infty , \,\  \hbox{where}\,\ \mu >0.
\end{equation}
We further remark from \cite[Proposition 3.1]{MV} that any solution $u$ of (\ref{eq:He}) satisfies the following Pohozaev-type identity
\begin{equation}\label{eq:pt}
 \frac{d-2}{2}\int_{\R^d}|\nabla u(x)|^2dx+\frac{d}{2}\int_{\R^d}|u(x)|^2dx=\frac{d-1}{2}\int_{\R ^d}\int_{\R ^d} \frac{u^2(x)u^2(y)}{|x-y|^2}dxdy.
\end{equation}
Following this identity, one can obtain that any solution $u(x)$ of (\ref{eq:He}) satisfies
\begin{equation}\label{eq:11}
\int_{\R ^d}|\nabla u(x)|^2dx=\int_{\R ^d}u^2(x)dx=\frac{1}{2}\int_{\R ^d} \frac{u^2(x)u^2(y)}{|x-y|^2}dxdy.
\end{equation}
Combining \eqref{eq:12} and \eqref{eq:11}, we know that all ground states $Q>0$ of (\ref{eq:He}) share the same $L^2$-norm, i.e.,
\begin{equation}\label{eq:N^*}
 N^*=\int_{\R^d}Q^2(x)dx,\ \forall \ Q\in G.
\end{equation}
We finally note from \cite{DMC}   the following Gagliardo-Nirenberg inequality
\begin{equation}\label{eq:GN}
\int_{\R ^d}\int_{\R ^d} \frac{u^2(x)u^2(y)}{|x-y|^2}dxdy\leq\frac{2}{N^*}\int_{\R ^d}|\nabla u(x)|^2dx\int_{\R ^d}u^2(x)dx,\,\ u\in H^1(\R^d),
\end{equation}
where the identity is attained at any ground state $Q=Q(|x|)$ of (\ref{eq:He}).
Denote $S_d>0$ the optimal constant of the following Sobolev  inequality
\begin{equation}\label{1:sob}
\Big(\int_{\R^d}|u|^{2^*}dx\Big)^{\frac{2}{2^*}}\leq S_d^{-1}\int_{\R^d}|\nabla u|^2dx,\,\ u\in H^1(\R^d),
\end{equation}
where $2^*=\frac{2d}{d-2}.$ Following above notations, our first main result of this paper is concerned with the existence and nonexistence of minimizers in terms of parameters $N$ and $\lambda$.

\begin{figure}
\centering
{\includegraphics[width = 8cm,height=5cm,clip]{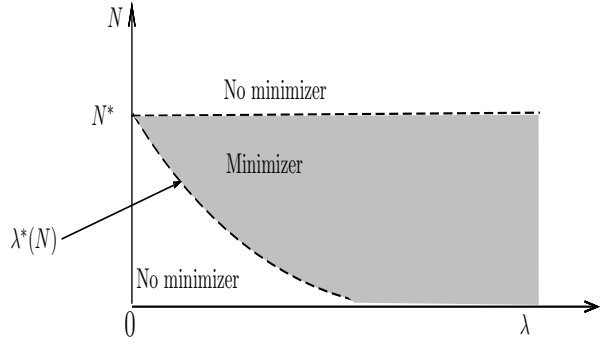}}
    \caption{\em  Existence  of minimizers for $e_\lambda (N)$ at different value $(\lam, N)$.}
 \end{figure}

\begin{thm}\label{thm1.1} Suppose $g(x)$ satisfies the assumption $\rm(M_1)$, and let $N^*:=\|Q\|^2_2>0$ be defined by (\ref{eq:N^*}). Then we have
\begin{enumerate}
\item [\rm(i).] If $N\geq N^*$, there is no minimizer of $e_\lambda(N)$ for any $\lambda>0$.
\item [\rm(ii).] If $0<N<N^*$, then there exists a positive constant $\lambda^*(N)$ satisfying
$$\frac{S_d(N^*-N)}{N^*\|1-g\|_{L^{\frac{d}{2}}}}\leq\lambda^*(N)\leq \frac{N^*-N}{\int_{\R^d}(1-g)Q^2dx},$$
such that
\begin{enumerate}
\item If $\lambda>\lambda^*(N)$, then there exists at least one minimizer for
$e_\lambda(N)$.
\item If  $0<\lambda<\lambda^*(N)$, then there is no minimizer for $e_\lambda(N)$.
\end{enumerate}
\end{enumerate}
\end{thm}

The existence and non-existence of minimizers in Theorem \ref{thm1.1} can be illustrated explicitly by Figure 1.  Theorem \ref{thm1.1} reveals the following interesting phenomenon: the depth $\lam$ of the steep potential well $g(x)$ must be large enough to keep the mass of minimizing sequences in a compact domain of $\R^d$, which ensures the existence of minimizers for $e_\lambda(N)$.
We shall prove that the finite constant $\lambda^*(N)$ in Theorem \ref{thm1.1} can be characterized by
\begin{equation}\label{def: boundary}
\lambda^*(N)=\underset{u\in H^1(\R^d),\|u\|^2_2=N}{\inf}\frac{\int_{\R^d}|\nabla u|^2dx-\frac{1}{2}\int_{\R^d}\int_{\R^d}\frac{u^2(x)u^2(y)}{|x-y|^2}dxdy}{\int_{\R^d}(1-g)u^2dx},
\end{equation}
where the assumption $1-g(x)\in L^{\frac{d}{2}}(\R^d)$ is to guarantee that $\lambda^*(N)>0$. Indeed, if $g(x)$ satisfies $\underset{|x|\rightarrow\infty}{\liminf }\big(1-g(x)\big)|x|^2=+\infty$, then one can get from (\ref{def: boundary}) by scaling that $\lambda^*(N)\equiv0$, which further implies the existence of minimizers for $e_\lambda(N)$ at any $0<N<N^*$ and $\lam>0$.

Since the energy functional $E_\lambda(\cdot)$ is even, any  minimizer  $u_\lam$ of $e_\lam (N)$ must be either $u_\lam\ge 0$ or $u_\lam\le 0$ in $\R^d$, and $u_\lam$ satisfies the following Euler-Lagrange equation
\begin{equation}\label{AA:euler}
-\bigtriangleup u_{\lambda} +\lambda g(x)u_{\lambda} =\mu_\lambda u_\lambda +\Big(\int_{\R^d}\frac{u^2_\lambda(y)}{|x-y|^2}dy\Big)u_\lambda \,\ \hbox{in}\,\ \R^d,
\end{equation}
where $\mu _\lam\in\R$ is the associated Lagrange multiplier. By the maximum principle, one can further derive from (\ref{AA:euler}) that $u_\lam$  must be either positive or negative. Without loss of generality, in the following we only consider positive  minimizers  $u_\lam>0$ of $e_\lam (N)$. We shall focus on the limit behavior of positive minimizers for $e_\lambda(N)$ as $\lambda\rightarrow\infty$, where $0<N<N^*$ is arbitrary, and $g(x)$ satisfies $(M_1)$ and
\begin{itemize}
\item [$\rm(M_2).$]
$g(x)\in C^{\alpha}_{loc}(\R^d)\,\, \hbox{for some}\,\,\alpha \in (0,1)$ and $
g(x)=|x|^p[1+o(1)]\,\,\,\hbox{as}\,\,\,|x|\rightarrow0$,\, where  \ $p>0$.
\end{itemize}
Under the above assumptions, our second main result of this paper can be stated as the following theorem.

\begin{thm}\label{thm1.2}
Suppose $g(x)$ satisfies $\rm(M_1)$ and $\rm(M_2)$ for $p>0$, and $N\in(0,N^*)$ is arbitrary. Let $u_{\lambda _k}$ be a positive minimizer of $e_{\lam_k}(N)$, where  $\lam _k\nearrow \infty$ as $k\to\infty$. Then there exists a subsequence of $\{\lam_k\}$ (still denoted by $\{\lam _k\}$) such that
\begin{equation}\label{def:18}
  \underset{ k\to\infty}{\lim} \lambda_k^{-\frac{d}{2(2+p)}} u_{\lambda_k}\Big(\lambda_k^{-\frac{1}{2+p}} x+x_{\lambda_k}\Big)=w_0(x)\,\ \hbox{in\, $H^1({\R^d})$},
\end{equation}
where $x_{\lambda_k}\in\R^d$ is a maximal point of $u_{\lambda _k}$, which satisfies $\lim_{k\to\infty}\lambda_k^{\frac{1}{2+p}}x_{\lambda_k}=0$, and  $w_0>0$ is a ground state of the following equation
\begin{equation}\label{def:333}
 -\triangle w_0 +|x|^pw_0 =\mu w_0 +\Big(\int_{\R^d}\frac{w^2_0(y)}{|x-y|^2}dy\Big)\,w_0 \,\ \hbox{in\, $\R^d$}, \  \ w_0\in H^1({\R^d}).
\end{equation}
Here $\mu =\mu(N)\in\R$ is determined by the constraint $\int_{\R^d} w^2_0=N$.
\end{thm}

The proof of Theorem \ref{thm1.2} shows that $w_0>0$ in (\ref{def:18}) is essentially a positive minimizer of the following minimization problem:
\begin{equation}\label{A:1.1}
  e_{\infty}(N)=\underset{\{u\in \mathcal{H}, \|u\|^2_2=N\}}{\inf} E_{\infty}(u),
\end{equation}
where the space $\mathcal{H}$ is defined by
\begin{equation}\label{A:1.2}
  \mathcal{H}=\Big\{u\in H^1(\R^d):\inte  |x|^pu^2dx<\infty\Big\},\ \ p>0,
\end{equation}
and the energy functional $E_{\infty}(u)$ is of the form
\begin{equation}\label{A:1.3}
  E_{\infty}(u)=\inte |\nabla u|^2dx+\inte |x|^pu^2dx-\frac{1}{2}\inte\inte\frac{u^2(x)u^2(y)}{|x-y|^2}dxdy,\ \ u\in   \mathcal{H}.
\end{equation}
The arguments of \cite{DLS, GS} give that $e_{\infty}(N)$ admits positive minimizers $w_0>0$ if and only if $N\in(0, N^*)$, where $N^*=\|Q\|^2_2$ and $Q>0$ is still a ground state of (\ref{eq:He}). As illustrated in Proposition \ref{prop:3.1}, any positive minimizer $w_0>0$ of $e_{\infty}(N)$ must be radially symmetric and strictly decreasing in $|x|$.  Moreover, $e_{\infty}(N)$ admits a unique positive minimizer, for the cases where either $N>0$ is small enough (cf. \cite[Theorem 1.1]{GZZ}), or $N$ is close enough to $N^*$ for $p\ge 2$ and $d=4$ (cf. \cite{Luo}).

In spite of above facts, we cannot obtain enough information, including the uniqueness and non-degeneracy, on positive minimizers of $e_{\infty}(N)$ for all cases of $0<N<N^*$ and $d\ge 3$. This leads to some new and challenging difficulties in the proof of Theorem \ref{thm1.2}, comparing with those appeared in \cite{GDL,DLS,GS,GWZZ,GZZ} and references therein. Actually, the existing methods in the above mentioned papers are only applicable to the case where $N<N^*$ is close sufficiently to $N^*$. To prove Theorem \ref{thm1.2} for the general case where $0<N<N^*$ is arbitrary, we shall establish Lemmas 3.4 and 3.5 on different types of estimates. Moreover, we need to seek for a different approach of deriving the lower bound of $e_\lam (N)$, for which we shall borrow some ideas from \cite{CKN}.

This paper is organized as follows. Section 2 is devoted to the proof of Theorem \ref{thm1.1} on the existence and nonexistence of  minimizers. In Section 3 we first address some a priori estimates of positive minimizers for $e_\lam (N)$ as $\lam\to\infty$, after which we shall complete the proof of Theorem \ref{thm1.2}.

\section{Existence of Minimizers}
In this section, we shall complete the proof of Theorem \ref{thm1.1} under the assumption $(M_1)$.
For convenience, we first define a new minimization problem
\begin {equation}\label{def:EA}
 \hat{e}_\lambda(N):=\underset{\{u\in H^1(\R^d),\|u\|^2_2=N\}}{\inf} \hat{E}_\lambda(u),\,\ d\ge 3,
 \end {equation}
where the energy functional $ \hat{E}_\lambda(u)$ satisfies for any $\lam >0$,
\begin{equation}
\hat{E}_\lambda(u):=\int_{\R ^d}|\nabla u(x)|^2dx-\lambda \int_{\R ^d}[1-g(x)]u^2(x)dx-\frac{1}{2} \int_{\R ^d}\int_{\R ^d} \frac{u^2(x)u^2(y)}{|x-y|^2}dxdy,
 \end{equation}
and $0\le g(x)\le 1$ satisfies the assumption $(M_1)$.
Since  $\hat{E}_\lambda(u)=E_\lambda(u)-\lambda N $, to prove Theorem \ref{thm1.1} it suffices to address the existence and nonexistence of minimizers for $\hat{e}_\lambda(N)$, instead of $e_\lambda(N)$. We start with the following lemma.

\begin{lem}\label{lem:existence1} Suppose $g(x)$ satisfies  the assumption $\rm(M_1)$ and let $N^*:=\|Q\|^2_2>0$ be defined by (\ref{eq:N^*}). Then
\begin{enumerate}
\item  If $N\geq N^*$, then there is no minimizer of (\ref{def:EA}) for any $\lambda>0$.
\item   If $0<N<N^*,\,\hbox{then}\,\, \hat{e}_\lambda(N)\leq 0$ for any $\lambda>0$.
\end{enumerate}
\end{lem}


\noindent{\bf Proof.}
\rm(i). We first consider the case where $N>N^*$ and $\lam >0$ are arbitrary. Set
\begin{equation}\label{equ:tf}
u_\theta(x)=\Big(\frac{N}{N^*}\theta^d\Big)^\frac{1}{2}Q(\theta x),\,\ \theta>0,
\end{equation}
where $Q>0$ is a ground state of (\ref{eq:He}).
By applying the Gagliardo-Nirenberg inequality (\ref{eq:GN}), we then have for any $N>N^*$ and $\lam >0$,
\begin{equation*}
\begin{split}
\hat{e}_\lambda(N) &\leq \hat{E}_\lambda (u_\theta(x))  \\
    &\leq \int_{\R^d}\frac{N}{N^*}\theta^d|\nabla Q(\theta x)|^2\theta^2 dx-\frac{N^2}{2{N^*}^2}\int_{\R^d}\int_{\R^d}\frac{\theta^{2d}Q^2(\theta x)Q^2(\theta y)}{|x-y|^2}dxdy\\
    &=N\Big(1-\frac{N}{N^*}\Big)\theta^2 \,\,        \rightarrow -\infty \ \ \mbox{as} \, \ \ \theta\rightarrow \infty,
\end{split}
\end{equation*}
due to (\ref{eq:11}), which thus implies that $\hat{e}_\lambda(N)=-\infty $ in this case. This shows the nonexistence of minimizers for $\hat{e}_\lambda(N)$, where $N>N^*$ and $\lam >0$.

Next, if $N=N^*$, we follow again the Gagliardo-Nirenberg inequality  (\ref{eq:GN}) to derive that for any $\lam >0$,
\begin{equation}\label{eq:Ean^*1}
   \hat{E}_\lambda(u) \geq -\lambda \int_{\R^d}(1-g)u^2dx \geq -\lambda N^*, \,\ i.e.,\,\ \hat{e}_\lambda(N^*)\geq -\lambda N^*.
\end{equation}
On the other hand, choose a cutoff function $0\leq \varphi \in C_0^\infty(\R^d)$ such that $\varphi(x)=1$ for $|x|\leq 1$ and $\varphi(x)=0$ for $|x|\geq 2$.
For $\tau>0$, let
\begin{equation}\label{def:trial}
u_\tau(x)=A_\tau\frac{\tau^{\frac{d}{2}}}{\| Q\|_2}\varphi(x) Q(\tau x),
\end{equation}
and $A_{\tau}>0$ is chosen so that $\int_{\R^d}u^2_\tau(x)dx=N^*$.
Using (\ref{eq:ed}), we have for some $C>0,$
\begin{equation*}
N^*\le A_\tau^2 \leq N^*+O(e^{-C\tau}) \,\ \ \text{as} \,\ \tau\to\infty.
\end{equation*}
By (\ref{eq:ed}), one can calculate that
\begin{equation}\label{test:del}
  \int_{\R^d}|\nabla u_\tau|^2dx \leq \frac{A_{\tau}^2 \tau^2}{\| Q\|_2^2}\int_{\R^d}|\nabla Q|^2dx +O(\tau^2e^{-C\tau}) \,\ \ \text{as}\,\ \tau\to\infty.
\end{equation}
As for the nonlocal term, we have
\begin{align}
      &    \int_{\R ^d}\int_{\R ^d} \frac{u^2_{\tau}(x)u^2_{\tau}(y)}{|x-y|^2}dxdy\nonumber\\
    = &   \frac{A_{\tau }^{4}\tau^2}{\|Q\|_2^{4}}
          \int_{\R^d}\int_{\R^d}\frac{\varphi^2(\frac{x}{\tau})\varphi^2(\frac{y}{\tau})Q^2(x)Q^2(y)}{|x-y|^{2}}dxdy\nonumber\\
    = &   \frac{A_{\tau }^{4}\tau^2}{\|Q\|_2^{4}}
\int_{\R^d}\int_{\R^d}\frac{Q^2(x)Q^2(y)}{|x-y|^2}dxdy
+O(\tau^2e^{-C\tau})\,\ \ \text{as}\,\ \tau\to\infty \label{cau:har}.
\end{align}
By (\ref{test:del}) and (\ref{cau:har}), we obtain from (\ref{eq:11}) that
\begin{equation}\label{eq:Ean^*2}
\begin{split}
\hat{e}_\lambda(N^*) \leq \hat{E}_\lambda(u_\tau(x))&\leq -\lambda \int_{\R^d}(1-g)u^2_{\tau}(x)dx+O(\tau^2e^{-C\tau}) \\
&=\frac{A_{\tau }^{2}}{\|Q\|_2^{2}}\lambda \int_{\R^d}g\big(\frac{x}{\tau}\big)\varphi^2\big(\frac{x}{\tau}\big)Q^2(x)dx-\lambda N^*+O(\tau^2e^{-C\tau}) \\
&=-\lambda N^*+o(1) \,\,\   \hbox{as} \,\,\  \tau\rightarrow \infty.
\end{split}
\end{equation}
We then conclude from (\ref{eq:Ean^*1}) and (\ref{eq:Ean^*2}) that $\hat{e}_\lambda(N^*)=-\lambda N^*$. If there were a minimizer $u$, we would  have
$$\int_{\R^d}g(x)u^2(x)dx=0=\inf_{x\in \R^d}g(x),$$ and $$\int_{\R ^d}|\nabla u(x)|^2dx=\frac{2}{N^*} \int_{\R ^d}\int_{\R ^d} \frac{u^2(x)u^2(y)}{|x-y|^2}dxdy,$$
which is a contradiction, since $u$ would have compact support by the first equality, while it has to be equal to $Q(x)$ (up to a scaling) in view of the second equality. Hence, we also have the nonexistence of minimizers for $\hat{e}_\lambda(N)$ in the case where $N=N^*$ and $\lam >0$. \vskip 0.1truein

\rm(ii). Using the same trial function  (\ref{equ:tf}), one can derive that
\begin{equation*}
\hat{e}_\lambda(N) \leq \hat{E}_\lambda (u_\theta(x))= \hat{E}_\lambda \big((\frac{N}{N^*}\theta^d)^\frac{1}{2}Q(\theta x)\big)\leq N\theta^2\Big(1-\frac{N}{N^*}\Big).
\end{equation*}
Following above estimates by taking $\theta\rightarrow 0,$ we thus conclude that $\hat{e}_\lambda(N) \leq 0$ holds for any $0<N<N^*$ and $\lam >0$, and the proof is therefore complete.
\qed

We next employ the celebrated concentration-compactness lemma (cf. \cite{Lions1}) to prove the following existence of minimizers.

\begin{lem}\label{thm:key1}
If $\hat{e}_\lambda(N)<0$ holds for some $0<N<N^*$ and $\lam >0$, then there exist minimizers for (\ref{def:EA}).
\end{lem}

\noindent{\bf Proof.}
For any $0<N<N^*$ and $\lam >0$, we have
\begin{equation}\label{eq:bbb}
  \begin{split}
    \hat{E}_\lambda(u) & =\int_{\R ^d}|\nabla u(x)|^2dx-\lambda \int_{\R ^d}[1-g(x)]u^2(x)dx-\frac{1}{2} \int_{\R ^d}\int_{\R ^d} \frac{u^2(x)u^2(y)}{|x-y|^2}dxdy\\
& \geq (1-\frac{N}{N^*})\int_{\R ^d}|\nabla u(x)|^2dx-\lambda N,\,\ \forall \,u\in H^1(\R^d).
  \end{split}
\end{equation}
This implies that any minimizing sequence $\{u_n\}$ of $\hat{e}_\lambda(N)$ is bounded uniformly in $H^1(\R^d)$ and $\|u_n\|^{2}_2= N$. It then follows from \cite[Lemma III.1]{Lions1} that there exists a subsequence $\{u_{n_k}\}$ of $\{u_n\}$ such that either the compactness or the dichotomy or the vanishing occurs for the subsequence $\{u_{n_k}\}$.

We first prove that the vanishing
\begin{equation}\label{eq:van}
\limsup_{k\rightarrow\infty,\,y\in \R^d}\int_{B_{R(y)}}{u^2_{n_k}}(x)dx=0,\,\,\, \forall\, R<\infty
\end{equation}
cannot occur. On the contrary, suppose that (\ref{eq:van}) is true. Then, we infer from  \cite[Lemma 1.21]{WM} that
\begin{equation}\label{eq:vanish}
{u_{n_k}}\overset{k}{\rightarrow} 0 \ \hbox{ in $L^q(\R^d)$, where $2<q<2^*$}.
\end{equation}
Following \eqref{eq:vanish}, since $g(x)$ satisfies $\rm(M_1)$, we obtain that
\begin{equation}\label{eq:zero}
  \inte [1-g(x)]u^2_{n_k}dx\leq \inte [1-g(x)]^{\frac{1}{2}}u^2_{n_k}dx\leq \|1-g\|_{\frac{d}{2}}^{\frac{1}{2}}\|u_{n_k}\|^2_{\frac{2d}{d-1}}\overset{k}{\rightarrow} 0.
\end{equation}
Similarly, by the Hardy-Littlewood-Sobolev inequality (c.f. \cite[Theorem 4.3]{EhL1}), we derive from (\ref{eq:vanish}) that
\begin{equation}\label{eq:hls}
\int_{\R ^d}\int_{\R ^d} \frac{u_{n_k}^2(x)u_{n_k}^2(y)}{|x-y|^2}dxdy\leq C\|u_{n_k}\|^4_{\frac{2d}{d-1}}\rightarrow0\ \ \hbox{as}\ \,k\rightarrow\infty.
\end{equation}
Using (\ref{eq:zero}) and (\ref{eq:hls}), we conclude that $\hat{e}_\lambda(N) =0$ in view of Lemma \ref{lem:existence1}(2),  which however contradicts the hypothesis that $\hat{e}_\lambda(N)<0$. This proves the claim (\ref{eq:van}).

To rule out the dichotomy, we now claim that for any $N\in (0,N^*)$,
 \begin{equation}\label{eq:sub}
   \hat{e}_\lambda(\theta\alpha)<\theta\, \hat{e}_\lambda(\alpha),\,\,\forall\, \alpha\in (0,N),\,\,\forall\,\theta\in \big(1,\frac{N}{\alpha}\big].
 \end{equation}
%
Indeed,   for any $\lambda>0$ and $0<N<N^*$, the estimate (\ref{eq:bbb}) implies that $\hat{e}_\lambda(\alpha)$ is bounded from below for any $\alpha\in (0,N)$. Let $\{v_n\}\subset H^1(\R^d)$ be a minimizing sequence of $\hat{e}_\lambda(\alpha)$, so that for any $n\in \mathbb{N}$, we have $ \|v_n\|^2_{2}=\alpha$ and $\lim\limits_{n\rightarrow\infty}\hat{E}_\lambda(v _n)=\hat{e}_\lambda(\alpha)$. It then follows from the estimate (\ref{eq:bbb}) that ${v_n}$ is bounded uniformly in $H^1(\R^d)$. Set $v_{\theta,n}:=\theta^{\frac{1}{2}}v_n(x)$ such that $\|v_{\theta,n}\|^2_{2}=\theta\alpha$. Since $\theta >1$, we obtain from (\ref{def:EA}) that
\begin{equation}\label{eq:bc}
\begin{split}
\hat{e}_\lambda(\theta\alpha)  \leq &\lim_{n\rightarrow\infty}\hat{E}_\lambda(v_{\theta,n}) \\
  =&\lim_{n\rightarrow\infty} \Big\{\theta\int_ {\R ^d}|\nabla v_n(x)|^2dx-\theta\lambda \int_{\R ^d}[1-g(x)]v_n^2(x)dx \\ &-\frac{1}{2}\theta^2 \int_{\R ^d}\int_{\R ^d}\frac{v_n^2(x)v_n^2(y)}{|x-y|^2}dxdy \Big\} \\
  =&\theta\hat{e}_\lambda(\alpha)+\frac{\theta(1-\theta)}{2} \lim_{n\rightarrow\infty}\int_{\R ^d}\int_{\R ^d}\frac{v_n^2(x)v_n^2(y)}{|x-y|^2}dxdy.
\end{split}
\end{equation}
If the minimizing sequence $\{v_n\}$ of $\hat{e}_\lambda(\alpha)$ satisfies
\begin{equation}\label{eq:minimi seq}
\liminf_{n\rightarrow\infty}\int_{\R ^d}\int_{\R ^d}\frac{v_n^2(x)v_n^2(y)}{|x-y|^2}dxdy\geq C_0>0,
\end{equation}
then the claim \eqref{eq:sub} follows in view of (\ref{eq:bc}). In order to prove \eqref{eq:minimi seq}, since $\hat {e}_{\lambda}(N)<0$, we obtain from the Gagliardo-Nirenberg inequality (\ref{eq:GN}) that
\begin{align*}
&\quad\lambda\liminf_{n\rightarrow\infty} \int_{\R ^d}[1-g(x)]v_n^2(x)dx\\
&=\liminf_{n\rightarrow\infty} \Big\{\int_ {\R ^d}|\nabla v_n(x)|^2dx-\frac{1}{2}\int_{\R ^d}\int_{\R ^d}\frac{v_n^2(x)v_n^2(y)}{|x-y|^2}dxdy\Big\}-\hat {e}_{\lambda}(N) \\
 & \geq \big(1-\frac{N}{N^*}\big)\liminf_{n\rightarrow\infty} \int_ {\R ^d}|\nabla v_n(x)|^2dx-\hat {e}_{\lambda}(N)>0.
\end{align*}
Following the above inequality, since $1-g(x)\in L^{\frac{d}{2}}(\R^d)$, we deduce that there exist constants $\varepsilon>0$ and large $R>0$ such that
\begin{equation}\label{eq:renon}
 \liminf_{n\rightarrow\infty} \int_{B_R(0)} v_n^2(x)\geq \varepsilon
\end{equation}
from which the estimate \eqref{eq:minimi seq} follows, and the claim \eqref{eq:sub} is hence established.

Once (\ref{eq:sub}) holds, we can follow \cite[Lemma II.1]{Lions1} to obtain that the following sub-additivity condition holds
\begin{equation}\label{eq:subadd}
  \hat{e}_\lambda(N)<\hat{e}_\lambda(\alpha)+\hat{e}_\lambda(N-\alpha).
\end{equation}
By contradiction, we now suppose that the  dichotomy occurs. It then follows from the arguments of \cite[Lemma III.1]{Lions1} that for $k>0$ large enough, there exist two sequences $\{u_{n_{k,1}}\}$ and $\{u_{n_{k,2}}\}$ such that $u_{n_k}$ satisfies
$$\|u_{n_k}-(u_{n_{k,1}}+u_{n_{k,2}})\|_q\rightarrow0 \ \,\hbox{as}\,\ k\rightarrow\infty\ \  \hbox{for\ $2\leq q<2^*$},$$
$$\Big|\int_{\R^d}u^2_{n_{k,1}}-\alpha\Big|\rightarrow0 \ \ \hbox{and}\,\ \Big|\int_{\R^d}u^2_{n_{k,2}}-(N-\alpha)\Big|\rightarrow0\,\ \hbox{as}\,\ k\rightarrow\infty,$$
$$\hbox{dist(Supp $u_{n_{k,1}}$, Supp $u_{n_{k,2}}$)}\rightarrow\infty\,\ \hbox{as}\,\  k\rightarrow\infty,$$
$$ \underset{k\rightarrow\infty}{\liminf}\int_{\R^d}\big\{|\nabla u_{n_k}|^2-|\nabla u_{n_{k,1}}|^2-|\nabla u_{n_{k,2}}|^2\big\}dx\geq0.$$
These further imply that
\begin{equation*}
\aligned
&~~\int_{\R ^d}\int_{\R ^d} \frac{u_{n_k}^2(x)u_{n_k}^2(y)}{|x-y|^2}dxdy\\
&=\int_{\R ^d}\int_{\R ^d}\frac{u^2_{n_{k,1}}(x)u^2_{n_{k,1}}(y)}{|x-y|^2}dxdy
+\int_{\R ^d}\int_{\R ^d}\frac{u^2_{n_{k,2}}(x)u^2_{n_{k,2}}(y)}{|x-y|^2}dxdy+o(1)
\endaligned
\end{equation*}
 as $k\to\infty$, and therefore,
\begin{equation*}
\begin{split}
 \hat{e}_\lambda(N)   =& \lim_{k\rightarrow\infty}\hat{E}_\lambda(u_{n_k}) \\
  \geq &\underset{k\rightarrow\infty}{\lim}\Big[\int_{\R^d}\big(|\nabla u_{n_{k,1}}|^2+|\nabla u_{n_{k,1}}|^2\big)dx
 -\lam \int_{\R ^d}[1-g(x)]\big(u^2_{n_{k,1}}+u^2_{n_{k,2}}\big)dx \\
 &-\frac{1}{2}\int_{\R ^d}\int_{\R^d}\frac{u^2_{n_{k,1}}(x)u^2_{n_{k,1}}(y)}{|x-y|^2}dxdy-\frac{1}{2}\int_{\R ^d}\int_{\R ^d}\frac{u^2_{n_{k,2}}(x)u^2_{n_{k,2}}(y)}{|x-y|^2}dxdy\Big]+o(1) \\
   \geq &\hat{e}_\lambda(\alpha)+\hat{e}_\lambda(N-\alpha),
\end{split}
\end{equation*}
which however contradicts (\ref{eq:subadd}). Therefore, the  dichotomy cannot occur, either.

By the concentration-compactness lemma  (cf. \cite{Lions1}), we now conclude that only the compactness occurs for the subsequence $\{u_{n_k}\}$. This implies that there exist a subsequence still denoted by $\{u_{n_k}\}$  and $\{y_k\}$ such that the sequence $\hat{u}_k(\cdot):=u_{n_k}(\cdot+y_k)$
satisfies
$$\hat u_{n_k}\rightharpoonup u_0 \,\,\,\hbox{ weakly in}\,\,\,H^1(\R^d),$$
and
$$\hat{u}_k\rightarrow u_0 \,\,\,\hbox{ strongly in}\,\,\,L^q(\R^d)\,\,\,\hbox{for}\,\,\, 2\leq q<2^*
$$
for some $u_0\in H^1(\R^d)$. So, we have
$$\lim_{k\rightarrow\infty}\int_{\R ^d}\int_{\R ^d} \frac{\hat{u}_k^2(x)\hat{u}_k^2(y)}{|x-y|^2}dxdy=\int_{\R ^d}\int_{\R ^d}\frac{u^2_0(x)u^2_0(y)}{|x-y|^2}dxdy.$$
By the weak lower semi-continuity, we also get that
$\liminf_{k\rightarrow\infty}\|\nabla \hat{u}_{k}\|^2_2\geq \|\nabla u_0\|^2_2,$  and so $u_0$ must be a minimizer of $\hat {e}_{\lambda}$, which further implies that $\lim_{k\rightarrow\infty}\|\nabla \hat{u}_{k}\|^2_2= \|\nabla u_0\|^2_2$.
Next, we claim that the translation $\{y_k\}$ is bounded uniformly for $k$. Actually, if the claim is false, we then have
\begin{align*}
\hat {e}_{\lambda }(N)&=\lim_{k\to\infty} \hat {E}_{\lambda }(u_{n_k})\\
&=\lim_{k\to\infty}\Big(\inte |\nabla u_{n_k}|^2dx-\lambda \inte[1-g(x)]u^2_{n_k}dx-\frac{1}{2}\inte\frac{u_{n_k}^2(x)u_{n_k}^2(y)}{|x-y|^2}dxdy\Big)\\
& =\lim_{k\to\infty}\Big(\inte |\nabla \hat{u}_{n_k}|^2dx-\lambda \inte[1-g(x+y_k)]\hat{u}^2_{n_k}dx
-\frac{1}{2}\inte\frac{\hat{u}_k^2(x)\hat{u}_k^2(y)}{|x-y|^2}dxdy\Big)\\
&=\inte |\nabla {u}_{0}|^2dx-\frac{1}{2}\inte\frac{{u}_0^2(x){u}_0^2(y)}{|x-y|^2}dxdy>\hat {E}_{\lambda }(u_0)=\hat {e}_{\lambda }(N),
\end{align*}
which is a contradiction. Once $\{y_k\}$ is bounded uniformly for $k$, taking a subsequence if necessary, we may assume that $\lim_{k\to\infty}y_k=y_0$, from which we obtain that
$$u_{n_k}\rightarrow u_0(\cdot-y_0)\ \ \hbox{ strongly in}\,\ L^q(\R^d)\ \ \hbox{as}\ \,k\rightarrow\infty, \,\,\ \hbox{where}\,\,\  2\leq q<2^*.$$
This implies that $u_0(\cdot-y_0)$ is therefore a minimizer of $\hat e_\lambda(N)$, and we are done. \qed

We next define $\lambda^*(N)$ by
\begin{equation}\label{eq:21}
  \lambda^*(N)=\underset{\{u\in H^1(\R^d),\,\|u\|^2_2=N\}}{\inf}F(u), \,\ 0<N<N^*,
\end{equation}
where $F(u)$ satisfies
\begin{equation}\label{eq:22}
F(u)=\frac{\int_{\R^d}|\nabla u|^2dx-\frac{1}{2}\int_{\R^d}\int_{\R^d}\frac{u^2(x)u^2(y)}{|x-y|^2}dxdy}{\int_{\R^d}(1-g)u^2dx}.
\end{equation}
The following lemma gives the existence and estimates of $\lambda^*(N)$.

\begin{lem}\label{lem:existence3}
For any $0<N<N^*$, let $\lambda^*(N)$ be defined by (\ref{eq:21}), where $g(x)$ satisfies $(M_1)$. Then $\lambda^*(N)$ satisfies \begin{equation}\label{2:est}0<\frac{S_d(N^*-N)}{N^*\|1-g\|_{\frac{d}{2}}}\leq\lambda^*(N)\leq \frac{N^*-N}{\int_{\R^d}(1-g)Q^2dx}<\infty,
\end{equation}
where $S_d>0$ is the optimal constant of the  Sobolev  inequality (\ref{1:sob}).
\end{lem}

\noindent{\bf Proof.}
By the definition of $\lambda^*(N)$, we reduce from (\ref{eq:11}) that
\begin{equation*}
 \lambda^*(N)\leq  F\big(\sqrt{\frac{N}{N^*}}Q\big)  =\frac{(N^*-N)\int_{\R^d}|\nabla Q|^2dx}{N^*\int_{\R^d}(1-g)Q^2dx}=\frac{N^*-N}{\int_{\R^d}(1-g)Q^2dx},
\end{equation*}
which implies the upper estimate of (\ref{2:est}).
As for the lower estimate of (\ref{2:est}), by the Gagliardo-Nirenberg inequality (\ref{eq:GN}), we derive that for any $u\in H^1(\R^d)$,
\begin{equation*}
  F(u)    \geq  \frac{(1-\frac{N}{N^*})\int_{\R^d}|\nabla u|^2dx}{\int_{\R^d}(1-g)u^2dx} \geq  \frac{(1-\frac{N}{N^*})\int_{\R^d(\R^d)}|\nabla u|^2dx}{\|1-g\|_{{\frac{d}{2}}}\|u\|^2_{{2^*}}}   \geq   \frac{S_d(N^*-N)}{N^*\|1-g\|_{\frac{d}{2}}},
\end{equation*}
where the Sobolev inequality (\ref{1:sob}) is used in the last inequality.
We thus obtain from above that the lower estimate of (\ref{2:est}) holds, and the lemma is therefore proved.
\qed
\vskip 0.1truein

\noindent{\bf Proof of Theorem \ref{thm1.1}.}
By Lemma \ref{lem:existence1}\rm(i), in the following we only need to prove Theorem \ref{thm1.1}\rm(ii).

Firstly, we consider the case where $0<N<N^*\, \hbox{and}\, \lambda>\lambda^*(N)$. By setting $\eta=\frac{\lambda-\lambda^*(N)}{2}>0$, there exists $u_{0}\in H^1(\R^d)$ satisfying $ \|u_{0}\|^2_2=N$ such that $$\lambda^*(N)\leq F(u_{0})<\lambda^*(N)+\frac{\eta}{2},$$
from which we infer that
\begin{equation}\label{eq:key2}
\begin{split}
\hat{e}_\lambda(N) & \leq\hat{E}_\lambda(u_0)<\big(\lambda^*(N)+\frac{\eta}{2}-\lambda\big)\int_{\R^d}(1-g)u^2_{0} \\
                   &=-\frac{3\eta}{2}\int_{\R^d}(1-g)u^2_{0}<0.
\end{split}
\end{equation}
Applying Theorem \ref{thm:key1}, we then conclude from (\ref{eq:key2}) that there exists a minimizer for (\ref{def:EA}) in this case, which therefore gives the existence of Theorem \ref{thm1.1}\rm(ii).

Next, we consider the case where $0<N<N^*\,\hbox{and}\,\, 0<\lambda<\lambda^*(N)$. By the definition of $\lambda^*(N)$, in this case we note that
 \begin{equation*}
\begin{split}
\hat{E}_\lambda(u) &\geq \lambda^*(N)\int_{\R^d}(1-g)u^2dx-\lambda\int_{\R^d}(1-g)u^2dx  \\
&=\big(\lambda^*(N)-\lambda\big)\int_{\R^d}(1-g)u^2dx\ge 0,
\end{split}
\end{equation*}
which implies that $\hat{e}_\lambda(N)\geq0$.
Since Lemma \ref{lem:existence1}\rm(ii) gives that $\hat{e}_\lambda(N)\leq0$, we have $\hat{e}_\lambda(N)=0$.
If (\ref{def:EA}) has a minimizer $u_0$, we then have
$$0=\hat{e}_\lambda(N)=\hat{E}_\lambda(u_0)\geq(\lambda^*(N)-\lambda)\int_{\R^d}(1-g)u^2_0dx\geq0,$$
which implies that $\int_{\R^d}(1-g)u^2_0dx=0$. However, this is impossible, since $1-g\ge 0$ and $g(x)\not =1$ in $\R^d$, and $\|u_0\|^2_2=N$. Therefore, we have the nonexistence of minimizers for the case where $0<N<N^*\,\hbox{and}\,\, 0<\lambda<\lambda^*(N)$. This establishes the nonexistence of Theorem \ref{thm1.1}\rm(ii), and we are done.
\qed

\section{Limit Behavior of Minimizers as $\lam \to\infty$}
This section is devoted to the proof of Theorem \ref{thm1.2} on the limit behavior of positive minimizers for $e_\lam (N)$ as $\lam \to\infty$, for the case where $0<N<N^*$, and $g(x)\in C^{\alpha}_{loc}(\R^d)$ satisfies $\rm(M_1)$ and $\rm(M_2)$ for $p>0$. In this case, we recall from Theorem \ref{thm1.1} that there exists at least one minimizer $u_\lambda$ for $e_\lam (N)$, i.e.,
\begin{equation*}
      e_\lambda(N)=\underset{\{u\in H^1(\R^d), \|u\|^2_2=N\}}{\inf} E_\lambda(u)=E_\lambda(u_\lambda).
\end{equation*}
We also introduce the following related minimization problem, which can be thought of as the ``limit" case of $e_\lam (N)$ for $\lam \to\infty$:
\begin{equation}\label{re3.2}
  e_{\infty}(N)=\underset{\{u\in \mathcal{H}, \|u\|^2_2=N\}}{\inf} E_{\infty}(u),
\end{equation}
where the space $\mathcal{H}$ is defined by \begin{equation}\label{A3}
  \mathcal{H}=\Big\{u\in H^1(\R^d):\inte  |x|^pu^2dx<\infty\Big\},\ \ p>0,
\end{equation}
and the energy functional $E_{\infty}(u)$ is of the form
\begin{equation}\label{re3.1}
  E_{\infty}(u)=\inte |\nabla u|^2dx+\inte |x|^pu^2dx-\frac{1}{2}\inte\inte\frac{u^2(x)u^2(y)}{|x-y|^2}dxdy,\ \ u\in   \mathcal{H}.
\end{equation}

Similar to \cite{DLS, GS}, one can prove that $e_{\infty}(N)$ admits positive minimizers if and only if $N\in(0, N^*)$, where $N^*=\|Q\|^2_2$ and $Q>0$ is a ground state of (\ref{eq:He}). We leave the detailed proof to the interested reader.
Moreover, we have the following qualitative properties of minimizers for $e_{\infty}(N)$, which play an important role in analyzing the limiting behaviour of minimizers $u_\lambda$ for $e_\lambda(N)$.


 \begin{prop}\label{prop:3.1}
 For any $N\in(0, N^*)$,  the following properties hold true:
\begin{enumerate}
  \item $\frac{e_{\infty}(N)}{N}\leq C(p)$, where $C(p)$ depends only on $p$.
  \item Any positive minimizer $U_N$ of $e_{\infty}(N)$ must be radially symmetric and strictly decreasing in $|x|$. Moreover, $U_N$ decays exponentially as $|x|\to\infty$.
\end{enumerate}
\end{prop}
\noindent{\bf Proof.} 1. By a transform $v= uN^{-\frac{1}{2}}$, for any $0<N<N^*$ the minimization problem $e_\lambda(N)$ is reduced equivalently to the following one:
\begin{equation}\label{A4}
  e_N=\underset{\{v\in \mathcal{H},\|v\|^2_2=1\}}{\inf} \mathcal{E}_N(v),
\end{equation}
where the energy functional $\mathcal{E}_{N}(v)$ is defined by
\begin{equation}\label{A5}
  \mathcal{E}_{N}(v)=\inte |\nabla v|^2dx+\inte |x|^pv^2dx-\frac{N}{2}\inte\inte\frac{v^2(x)v^2(y)}{|x-y|^2}dxdy,\ \ p>0.
\end{equation}
It is clear that $e_N=\frac{e_{\infty}(N)}{N}$. By the arguments of \cite{DLS,GS}, one can reduce that there exists a constant $C(p)>0$ depending only on $p$ such that $e_{N}\leq C(p)$, and hence $\frac{e_{\infty}(N)}{N}\leq C(p)$ for any $0<N<N^*$.

2. Using the symmetric-decreasing rearrangement, similar to \cite{CWZ,EHL2} one can obtain that any positive minimizer $U_N(x)$ of $e_{\infty}(N)$ must be radially symmetric and decreasing in $|x|$. Further, similar to \cite{GS,GWZZ,GZZ}, the comparison principle yields that $U_N(x)$ decays exponentially as $|x|\to\infty$.

We finally prove that $U_N(x)$ strictly decreases in $|x|$. On the contrary, suppose that there exist  two points $y_1\in\R^d$ and $y_2\in\R^d$ such that $|y_1|\neq |y_2|$ and $U_N(y_1)=U_N(y_2)$. Without loss of generality, we may assume $0\le a=|y_1|<b=|y_2|$. Since $U_N(x)>0$ decreases in $|x|$, we have $U_N(x)\equiv const.>0$ in the annual domain $\bar {\Omega}:=\bar {B}_b(0)/B_a(0)$. Hence, we derive from (\ref{def:333}) that $-\triangle U_N(x)\equiv 0$ in $\bar {\Omega}$ and
\begin{equation}\label{A2*}
|x|^p U_N(x)=\mu  U_N(x)+\Big(\int_{\R^d}\frac{ U_N^2(y)}{|x-y|^2}dy\Big)U_N(x)\quad\hbox{in\, $\bar {\Omega}$.}
\end{equation}
Applying \eqref{A2*}, we have
\begin{equation}\label{A3*}
  f(x)=\int_{\R^d}\frac{U_N^2(y)}{|x-y|^2}dy-|x|^p+\mu\equiv 0\quad\hbox{in\, $\bar {\Omega}$,}
\end{equation}
which thus gives that
\begin{equation}\label{A4*}
f(P_1)= f(P_2),\quad\hbox{where $P_1=(a,0,\cdots,0), P_2=(b,0,\cdots,0)$.}
\end{equation}
On the other hand, for any point $y=(y_1,y_2,\cdots,y_d)\in\R^d$, set $P=(\frac{a+b}{2},0,\cdots,0)$
and define $y_P=(a+b-y_1,y_2,\cdots,y_d)$ such that
\[
|P_1-y|^2\leq|P_2-{y}|^2\ \, \mbox{and}\,\  |y|<|y_P|\ \ \mbox{in}\,\ \big\{y\in \R^d:\, y_1\leq \frac{a+b}{2}\big\}.
\]
By direct calculations, we then deduce that
\begin{equation}\label{A6*}
  \begin{split}
&\quad  f(P_1)-f(P_2)\\
&=\inte \Big[\frac{U_{N}^2(y)}{|P_1-y|^2}-\frac{U_{N}^2(y)}{|P_2-y|^2}\Big]dy+(b^q-a^q)\\
&=\int_{\{y\in \R^d: y_1\leq\frac{a+b}{2}\}} U_{N}^2(y)\Big[\frac{1}{|P_1-y|^2}-\frac{1}{|P_2-y|^2}\Big]dy\\
&~+\int_{\{y\in \R^d: y_1>\frac{a+b}{2}\}} U_{N}^2(y)\Big[\frac{1}{|P_1-y|^2}-\frac{1}{|P_2-y|^2}\Big]dy+(b^q-a^q)\\
&=\int_{\{y\in \R^d: y_1\leq\frac{a+b}{2}\}} U_{N}^2(y)\Big[\frac{1}{|P_1-y|^2}-\frac{1}{|P_2-y|^2}\Big]dy\\
&~+\int_{\{y\in \R^d: y_1<\frac{a+b}{2}\}} U_{N}^2(y_P)\Big[\frac{1}{|P_2-y|^2}-\frac{1}{|P_1-y|^2}\Big]dy+(b^q-a^q)\\
&=\int_{\{y\in \R^d: y_1\leq\frac{a+b}{2}\}} \Big[U_{N}^2(y)-U_{N}^2(y_P)\Big]\Big[\frac{1}{|P_1-y|^2}-\frac{1}{|P_2-y|^2}\Big]dy+(b^q-a^q)>0,\\
  \end{split}
\end{equation}
since $U_N(x)$ decreases to zero in $|x|$. We thus conclude from \eqref{A6*} that $f(P_1)>f(P_2)$, which however contradicts \eqref{A4*}.
Therefore, $U_N(x)$ strictly decreases in $|x|$, which completes the proof.
\qed


We next apply the above qualitative properties of  $e_{\infty}(N)$ to addressing the following refined energy estimates.


\begin{lem}\label{lem:energy estimate}
Suppose $g(x)$ satisfies $\rm(M_1)$ and $\rm(M_2)$ for $p>0$. Then for any $N\in(0, N^*)$, we have
\begin{equation}\label{eq:energy estimate}
B(p)\Big(1-\frac{N}{N^*}\Big)^{\frac{p}{2+p}}N\lambda^\frac{2}{2+p}\leq e_\lambda(N)\leq \lambda^\frac{2}{2+p}\big[e_{\infty}(N)+o (1)\big]\ \, \hbox{as}\,\ \lambda\rightarrow\infty,
\end{equation}
where the constant $B(p)>0$ depends only on $p$.
\end{lem}

\noindent{\bf Proof.} By Proposition \ref{prop:3.1}, let $U_N>0$ be a positive minimizer of $e_{\infty}(N)$ for any $N\in(0, N^*)$, where the problem $e_{\infty}(N)$ is defined in \eqref{re3.2}. Take the trial function $\hat{U}_N(x)= \lambda^{\frac{d}{2(2+p)}}U_N(\lambda^{\frac{1}{2+p}}x)$.
We then have
\begin{equation}\label{eq:31}
\begin{split}
&\quad e_\lambda(N)\leq E_\lambda(\hat{U}_N) \\
&=\lambda^{\frac{2}{2+p}}\Big[\inte |\nabla {U}_N(x)|^2dx+\lambda^{\frac{p}{2+p}}\inte g\Big(\frac{x}{\lambda^{\frac{1}{2+p}}}\Big)U^2_N(x)dx\\
&\quad\quad \qquad -\frac{1}{2}\inte\inte\frac{{U}^2_N(x){U}^2_N(y)}{|x-y|^2}dxdy\Big]\\
&=\lambda^{\frac{2}{2+p}}\Big(e_{\infty}(N)+\inte \Big[\lambda^{\frac{p}{2+p}}g\Big(\frac{x}{\lambda^{\frac{1}{2+p}}}\Big)-|x|^p\Big]U^2_N(x)dx\Big)\ \, \hbox{as}\,\ \lambda\rightarrow\infty.
\end{split}
\end{equation}
Since $g(x)$ satisfies $\rm(M_1)$ and $\rm(M_2)$ for $p>0$, we deduce that there exists a large constant $C>0$ such that
$$g(x)\leq C|x|^p  \ \ \mbox{and}\,\ \lambda^{\frac{p}{2+p}}g\Big(\frac{x}{\lambda^{\frac{1}{2+p}}}\Big)\leq C|x|^p\ \ \mbox{in}\,\ \R^d.$$
Since Proposition \ref{prop:3.1} gives that $U_N(x)$ decays exponentially as $|x|\to\infty$, we have
\[
\int_{|x|> \lambda^{\frac{1}{2(2+p)}}} \Big[\lambda^{\frac{p}{2+p}}g\Big(\frac{x}{\lambda^{\frac{1}{2+p}}}\Big)-|x|^p\Big]U^2_N(x)dx=o(1)\,\ \hbox{as}\,\   \lambda\rightarrow\infty.
\]
We thus calculate that
\begin{align}\label{eq:gx}
&\quad \inte \Big[\lambda^{\frac{p}{2+p}}g\Big(\frac{x}{\lambda^{\frac{1}{2+p}}}\Big)-|x|^p\Big]U^2_N(x)dx\nonumber\\
&=\int_{|x|\leq \lambda^{\frac{1}{2(2+p)}}} \Big[\lambda^{\frac{p}{2+p}}g\Big(\frac{x}{\lambda^{\frac{1}{2+p}}}\Big)-|x|^p\Big]U^2_N(x)dx\\
&\quad+\int_{|x|> \lambda^{\frac{1}{2(2+p)}}} \Big[\lambda^{\frac{p}{2+p}}g\Big(\frac{x}{\lambda^{\frac{1}{2+p}}}\Big)-|x|^p\Big]U^2_N(x)dx\nonumber \\
&=\big(1+o (1)-1\big)\int_{|x|\leq \lambda^{\frac{1}{2(2+p)}}}|x|^pU^2_N(x)dx+o (1)=o(1)\,\ \hbox{as}\,\   \lambda\rightarrow\infty,\nonumber
\end{align}
from which we obtain the upper energy estimate of (\ref{eq:energy estimate}).


Next, we shall address the lower estimate of (\ref{eq:energy estimate}). By the Gagliardo-Nirenberg inequality (\ref{eq:GN}), we deduce from the upper energy estimate of (\ref{eq:energy estimate}) that
$$\lambda\int_{\R^d}g(x)u^2_\lambda dx\leq e_\lambda(N)\leq \big[e_{\infty}(N)+o (1)\big]\lambda^\frac{2}{2+p} \,\ \hbox{as}\,\   \lambda\rightarrow\infty,$$
which implies that
\begin{equation}\label{eq:38}
 \int_{\R^d}g(x)u^2_\lambda dx\rightarrow 0 \,\ \hbox{as}\,\   \lambda\rightarrow\infty.
\end{equation}
On the other hand, under the assumption $\rm(M_2)$ for $p>0$, let $\varepsilon_0>0$ be   small enough that $g(x)\geq \frac{1}{2}|x|^p$ holds for any $|x|\leq 2\varepsilon_0$. Because $0$ is the unique global minimum point of $g(x)$ satisfying $g(0)=0$, we obtain from (\ref{eq:38}) that
\[\int_{\{|x|> \varepsilon_0\}}u^2_\lambda dx\rightarrow 0\ \ \hbox{as}\,\  \lambda\rightarrow\infty,\]
and
\[\int_{\{|x|\leq  \varepsilon_0\}}u^2_\lambda dx=\inte{u}^2_\lambda-\int_{\{|x|> \varepsilon_0\}}{u}^2_\lambda dx\geq\frac{N}{4}\,\ \hbox{as}\,\   \lambda\rightarrow\infty.  \]
Define now $\hat{u}_\lambda=\varphi u_\lambda$, where $0\le \varphi \in C_0^\infty(\R^d) $ is a smooth cut-off function satisfying
\begin{equation}\label{eq:37}
\left\{
  \arraycolsep=1.5pt \begin{array}{ll}
     & \hbox{$\varphi=1\ \,\hbox{for}\,\ |x|\leq\varepsilon_0$;} \\[2mm]
    & \hbox{$\displaystyle 0<\varphi<1 \ \,\hbox{for}\,\  \varepsilon_0<|x|< 2\varepsilon_0$;}\\[2mm]
    & \hbox{$\varphi=0\ \,\hbox{for}\,\   |x|\geq 2\varepsilon_0$;} \,\ \hbox{$|\nabla \displaystyle \varphi|\leq\frac{M}{\varepsilon_0}\ \,\hbox{in}\,\  \R^d$.}
  \end{array}
\right.
\end{equation}
One can check that $\hat{u}_\lambda$ satisfies
\begin{equation}\label{eq:39}
\left\{  \arraycolsep=1.5pt
  \begin{array}{ll}
     & \hbox{$0\leq \hat{u}_\lambda(x)\leq u_\lambda(x)$ in $ \R^d$,} \\[2mm]
     & \hbox{$\displaystyle \frac{N}{4}\leq \|\hat {u}_\lambda\|^2_2\leq \|u_\lambda\|^2_2$,} \\[2mm]
    & \displaystyle\hbox{$\displaystyle \|\nabla\hat{u}_\lambda\|^2_2\leq2\|\nabla u_\lambda\|^2_2+\frac{2M^2}{\varepsilon^2_0}N$.}
  \end{array}
\right.
\end{equation}
By the Caffarelli-Kohn-Nirenberg inequality (cf. \cite{CKN}), we then reduce from (\ref{eq:39}) that there exists a constant $A(p)>0$ depending only on $p$ such that
\begin{equation}\label{eq:310}
\begin{split}
 \|\nabla\hat{u}_\lambda \|^{\frac{p}{2+p}}_2\|g^{\frac{1}{2}}(x)\hat{u}_\lambda\|^{\frac{2}{2+p}}_2 & \geq \big(\frac{1}{2}\big)^{\frac{1}{2+p}}\|\nabla\hat{u}_\lambda\|^{\frac{p}{2+p}}_2\| |x|^{\frac{p}{2}}\hat{u}_\lambda\|^{\frac{2}{2+p}}_2\\                                                                                               &\geq \big(\frac{1}{2}\big)^{\frac{1}{2+p}}A(p)\|\hat{u}_\lambda\|_2 \geq \displaystyle \big(\frac{1}{2}\big)^{1+\frac{1}{2+p}}A(p)N^{\frac{1}{2}}.
\end{split}
\end{equation}
Since it yields from (\ref{eq:38}) that
\[\inte g(x)\hat{u}^2_\lambda(x)\leq \inte g(x)u^2_\lambda(x)\to 0\,\ \hbox{as}\,\ \lambda\to\infty,\]
we obtain from \eqref{eq:310} that
\begin{equation}\label{re:310}
  \|\nabla\hat{u}_\lambda \|_2\to \infty\,\ \hbox{as}\,\ \lambda\to\infty.
\end{equation}
By \eqref{eq:39} and \eqref{re:310}, we then have
\begin{equation}\label{re:311}
\|\nabla u_\lambda \|^2_2\geq \frac{1}{2}\|\nabla\hat{u}_\lambda \|^2_2-\frac{M^2N}{\varepsilon^2_0}\geq \frac{1}{4}\|\nabla\hat{u}_\lambda \|^2_2\ \ \hbox{as}\,\ \lambda\to\infty.
\end{equation}
Following \eqref{eq:39}, (\ref{eq:310}) and (\ref{re:311}), we thus conclude that
\begin{equation}\label{eq:311}
\begin{split}
\|\nabla u_\lambda\|^{\frac{p}{2+p}}_2\|g^{\frac{1}{2}}(x)u_\lambda\|^{\frac{2}{2+p}}_2&\geq \Big(\frac{1}{4}\|\nabla\hat{u}_\lambda \|^2_2\Big)^{\frac{p}{2(2+p)}}\|g^{\frac{1}{2}}(x)\hat{u}_\lambda \|^{\frac{2}{2+p}}_2 \\
& \geq \Big(\frac{1}{2}\Big)^{1+\frac{p+1}{2+p}}A(p)N^{\frac{1}{2}}\ \ \hbox{as}\,\  \lambda\rightarrow\infty.
\end{split}
\end{equation}
By the Gagliardo-Nirenberg inequality (\ref{eq:GN}), together with Young's inequality, we now reduce from (\ref{eq:311}) that
\begin{equation}\label{eq:312}
 \begin{split}
   e_\lambda(N)&\geq \Big(1-\frac{N}{N^*}\Big)\|\nabla u_\lambda\|^2_2+\lambda \|g^{\frac{1}{2}}(x)u_\lambda\|^2_2 \\
&\geq  (2+p)(4p^p)^{-\frac{1}{2+p}}\Big(1-\frac{N}{N^*}\Big)^{\frac{p}{2+p}}\lambda^\frac{2}{2+p}\|\nabla u_\lambda\|^\frac{2p}{2+p}_2\|g^{\frac{1}{2}}(x)u_\lambda\|^\frac{4}{2+p}_2 \\
&\geq \frac{1}{16}(2+p)p^{-\frac{p}{2+p}}A^2(p)N\Big(1-\frac{N}{N^*}\Big)^{\frac{p}{2+p}}\lambda^{\frac{2}{2+p}}\\
&:=B(p)\Big(1-\frac{N}{N^*}\Big)^{\frac{p}{2+p}}N\lambda^{\frac{2}{2+p}}\ \ \hbox{as}\,\  \lambda\rightarrow\infty,
 \end{split}
\end{equation}
which then implies that the lower bound estimate of  (\ref{eq:energy estimate}) also holds.
\qed

Using (\ref{eq:energy estimate}) and (\ref{eq:311}), one can establish immediately the following corollary.

\begin{cor}\label{re: 313}
Suppose $g(x)$ satisfies $\rm(M_1)$ and $\rm(M_2)$ for $p>0$, and let $u_\lambda>0$ be a positive minimizer of $e_\lam (N)$. Then for any $N\in(0, N^*)$, we have
\begin{equation}\label{eq:313}
\left\{\arraycolsep=1.5pt
  \begin{array}{ll}
    &\hbox{$m_1(N,p)\lambda^\frac{2}{2+p}\leq \|\nabla u_\lambda \|^2_2\leq m_2(N,p)\lambda^\frac{2}{2+p}$ \,    as\,  $\lambda\rightarrow\infty$,} \\ [2mm]
    &\hbox{$M_1(N,p)\lambda^\frac{2}{2+p}\leq\lambda \|g^{\frac{1}{2}}(x)u_\lambda \|^2_2\leq M_2(N,p)\lambda^\frac{2}{2+p}$\,\   as\,  $\lambda\rightarrow\infty$,}
  \end{array}
\right.
\end{equation}
where $m_1(N,p)$, $m_2(N,p)$, $M_1(N,p)$ and $M_2(N,p)$ are positive constants depending only on $N$ and $p$.
\end{cor}

\noindent{\bf Proof.} Indeed, the upper bound estimates of (\ref{eq:313}) follow directly from Lemma \ref{lem:energy estimate} and  the Gagliardo-Nirenberg inequality (\ref{eq:GN}) as well.  As for the lower  bound estimates of (\ref{eq:313}), by contradiction we first suppose that $m_1(N, p)=0$, i.e., suppose  there exists a sequence $\{\lambda_n\}$, where $\lambda_n\to\infty$ as $n\to\infty$, such that
\begin{equation}\label{eq:314}
  \lim_{\lambda_n\rightarrow\infty}\frac{\|\nabla u_{\lambda_n}\|^2_2}{\lambda^{\frac{2}{2+p}}_n}=0.
\end{equation}
Following (\ref{eq:311}) and (\ref{eq:314}), we then obtain that
\begin{equation*}
\lim_{\lambda_n\rightarrow\infty}\frac{\lambda_n\|g^{\frac{1}{2}}(x)u_{\lambda_n}\|^2_2}{\lambda^\frac{2}{2+p}_n}=
\lim_{\lambda_n\rightarrow\infty}\|g^{\frac{1}{2}}(x)u_{\lambda_n}\|^2_2 \lambda^{\frac{p}{2+p}}_n=\infty,
\end{equation*}
which however contradicts   (\ref{eq:energy estimate}). Therefore, we have $m_1>0$. Similarly, one can also obtain that $M_1>0$, and therefore (\ref{eq:313}) is proved.
\qed

\vskip 0.1truein

Since $u_\lambda>0$ is a positive minimizer of $E_\lambda(u)$, it satisfies the following Euler-Lagrange equation
\begin{equation}\label{eq:euler}
-\bigtriangleup u_{\lambda} +\lambda g(x)u_{\lambda} =\mu_\lambda u_\lambda +\Big(\int_{\R^d}\frac{u^2_\lambda(y)}{|x-y|^2}dy\Big)u_\lambda \,\ \hbox{in}\,\ \R^d,
\end{equation}
where $\mu _\lam\in\R$ is the associated Lagrange multiplier satisfying
\begin{equation}\label{re:lagrange}
\mu_\lambda=\frac{1}{N}\Big[e_\lambda(N)-\frac{1}{2} \int_{\R ^d}\int_{\R ^d} \frac{u^2_\lambda(x)u^2_\lambda(y)}{|x-y|^2}dxdy\Big].
\end{equation}
We next establish the following estimates.

\begin{lem}\label{prop u} Suppose $g(x)$ satisfies $\rm(M_1)$ and $\rm(M_2)$ for $p>0$, and let $u_\lambda>0$ be a positive minimizer of $e_\lam (N)$. Then we have
\begin{enumerate}
\item  There exists a constant $C(p)> 0$, depending only on $p>0$, such that the Lagrange multiplier $\mu_\lambda$ of (\ref{re:lagrange}) satisfies
\begin{equation}\label{prop u1}
  \mu_{\lambda}\leq C(p)\,\lambda^{\frac{2}{2+p}}\,\ \hbox{as}\,\ \lambda\to\infty.
\end{equation}

\item  $u_{\lambda}$ has at least one maximal point $x_\lambda$, which satisfies
\begin{equation}\label{prop u2}
  u_{\lambda}(x_\lambda)\geq C_0(N,p)\,\lambda^{\frac{d}{2(2+p)}}\,\ \hbox{as}\,\ \lambda\to\infty
\end{equation}
for some constant $C_0(N,p)> 0$.
\end{enumerate}
\end{lem}

\noindent{\bf Proof.} (1). By Proposition \ref{prop:3.1} and Lemma \ref{lem:energy estimate}, we deduce from \eqref{re:lagrange} that
\begin{equation*}
\begin{split}
\mu_{\lambda} & =\frac{1}{N}\Big[e_\lambda(N)-\frac{1}{2} \int_{\R ^d}\int_{\R ^d} \frac{u^2_\lambda(x)u^2_\lambda(y)}{|x-y|^2}dxdy\Big]\\
              &\leq \frac{1}{N}e_\lambda(N)\leq \frac{1}{N}[e_{\infty}(N)+o (1)]\lambda^{\frac{2}{2+p}}\leq C(p)\lambda^{\frac{2}{2+p}}\,\ \hbox{as}\,\ \lambda\to\infty,
\end{split}
\end{equation*}
which then gives the estimate (\ref{prop u1}).

(2). Following \eqref{prop u1} yields that $\underset{\lambda\rightarrow\infty}{\liminf}[\lambda g(x)-\mu_\lam ]>0$ for $|x|>R$, where $R>0$ is large enough. Together with  (\ref{eq:euler}), we then have
\begin{equation}\label{re: 316}
  -\triangle u_{\lambda}-\Big(\int_{\R ^d}\frac{u^2_\lambda(y)}{|x-y|^2}dy\Big)u_\lambda\leq 0,\,\ \hbox{where $|x|\geq R>0$.}
\end{equation}
The Hardy-Littlewood-Sobolev inequality  (c.f. \cite[Theorem 4.3]{EhL1}) gives that $\inte\frac{u^2_\lambda(y)}{|x-y|^2}dy \in L^d(\R^d)$. Following the De Giorgi-Nash-Moser theory (c.f. \cite[Theorem 4.1]{HL}), we thus derive from (\ref{re: 316}) that for any $\xi\in B_R^c(0)$,
\begin{equation*}
  \max_{x\in B_1(\xi)}u_{\lambda}(x)\leq C\Big(\int_{B_2(\xi)}u^2_{\lambda}(x)dx\Big)^{\frac{1}{2}}.
\end{equation*}
The above inequality implies that $u_{\lambda}(x)$ tends to zero near infinity, and thus each $u_\lambda$ has at least one maximum point which we denote by $x_\lambda$.

We finally prove \eqref{prop u2} as follows. For convenience, denote $\varepsilon_{\lambda}=\lambda^{-\frac{1}{2+p}}>0$. By contradiction  suppose the estimate \eqref{prop u2} is false. Then for any small constant $\varepsilon>0$, we have
\begin{equation}\label{re:317}
   u_{\lambda}(x_\lambda)\leq \varepsilon\varepsilon^{-\frac{d}{2}}_{\lambda}\,\ \hbox{in\  $\R^d$.}
\end{equation}
On the other hand, for any large $l>0$, it follows from Corollary \ref{re: 313} that
\begin{equation}\label{re:341}
\int_{\{g(x)\geq l\varepsilon_{\lambda}^p\}}u^2_\lambda(x)dx\leq
\frac{1}{l\varepsilon_{\lambda}^p}\int_{\R^d}g(x)u^2_\lambda(x)dx\leq\frac{M_2}{l}.
\end{equation}
We thus deduce from \eqref{re:317} and \eqref{re:341} that for large $l>0,$
\begin{equation}\label{re:342}
\begin{split}
\frac{N}{2}\geq\frac{M_2}{l}&\geq\int_{\{g(x)\geq l\varepsilon_{\lambda}^p\}}u^2_\lambda(x)dx=N-\int_{\{g(x)\leq l\varepsilon_{\lambda}^p\}}u^2_\lambda(x)dx \\
             &\geq N-\int_{\{|x|\leq 2l^{\frac{1}{p}}\varepsilon_{\lambda}\}}u^2_\lambda(x)dx\geq N-\int_{\{|x|\leq 2l^{\frac{1}{p}}\varepsilon_{\lambda}\}}\varepsilon^2\varepsilon^{-d}_{\lambda}dx\\
             &\geq N-C\varepsilon^2> \frac{N}{2}\,\ \hbox{as}\,\ \lambda\to\infty.
\end{split}
\end{equation}
This is a contradiction, and hence \eqref{prop u2} holds true.  The proof is therefore complete.
\qed
\vskip 0.1truein

In view of the estimate (\ref{eq:energy estimate}), we now define
the $L^2(\R^d)$-normalized function
\begin{equation}\label{def:316}
  w_\lambda(x):=\varepsilon_{\lambda}^{\frac{d}{2}} u_\lambda(\varepsilon_\lambda x+x_\lambda), \,\ \varepsilon_\lambda:=\lambda^{-\frac{1}{2+p}}>0,
\end{equation}
where  $x_{\lambda}\in\R^d$ is a maximal point of $u_{\lambda}$.
We then deduce from (\ref{eq:313}) that
\begin{equation}\label{def:317}
 \|\nabla w_\lam \|^2_2\leq C(N,p),\,\,\,\|w_\lambda\|^2_2=N,
\end{equation}
and note from (\ref{eq:euler}) that $w_\lambda>0$ satisfies the following Euler-Lagrange equation
\begin{equation}\label{eq:Euler}
-\bigtriangleup w_{\lambda}(x)+\varepsilon_{\lambda}^2\lambda\, g(\varepsilon_{\lambda} x+x_\lambda)w_{\lambda}(x)=\varepsilon_{\lambda}^2\mu_{\lambda} w_\lambda(x)+\Big(\int_{\R^d}\frac{w^2_\lambda(y)}{|x-y|^2}dy\Big)w_\lambda(x) \,\ \hbox{in}\,\ \R^d.
\end{equation}

\begin{lem}\label{exp decaey}
Suppose $g(x)$ satisfies $\rm(M_1)$ and $\rm(M_2)$ for $p>0$, and let $u_\lambda>0$ be a positive minimizer of $e_\lam (N)$. Then we have
\begin{enumerate}
  \item  The maximal point $x_{\lambda}$ of $u_{\lambda}$ satisfies
\begin{equation}\label{point extimae}
  |x_{\lambda}|\leq C(N,p)\varepsilon_{\lambda} \ \, \hbox{as \  $\lambda\to\infty$},
\end{equation}
where $C(N,p)>0$ is independent of $\lambda$.

 \item  There exist sufficiently large $R=R(N, p)>0$ and $\lambda_0>0$ such that $w_{\lambda}$ satisfies
\begin{equation}\label{re exp}
w_\lambda(x)\leq C(N,p)e^{-|x|}\, \ \hbox{in}\,\ B^c_R(0),\,\  \mbox{if}\,\ \lambda>\lambda_0.
\end{equation}

\item  There exist a subsequence  $\{w_k\}$ of $\{w_{\lambda_k}\}$, where $\lambda_k\to\infty$ as $k\to\infty$, and $0\le w_0\in H^1(\R^d)$ such that
\begin{equation}\label{re conver}
w_k\rightharpoonup w_0 \,\ \hbox{in}\,\  H^1(\R^d) \ \, \hbox{as \  $k\to\infty$},
\end{equation}
and
\begin{equation}\label{re conver*}
w_k\rightarrow w_0\,\  \hbox{strongly in}\,\  L^q(\R^d) \ \, \hbox{as \  $k\to\infty$},\,\ \hbox{where\, $2\leq q<2^*$}.
\end{equation}

\end{enumerate}
\end{lem}

\noindent{\bf Proof.} (1). By contradiction, suppose \eqref{point extimae} is false, which implies that \begin{equation}\label{re:32T}\lim_{\lambda\to\infty} |x_{\lambda}|\varepsilon^{-1}_{\lambda}\to \infty,\ \ \mbox{where}\,\ \varepsilon_\lambda:=\lambda^{-\frac{1}{2+p}}>0.
\end{equation}
Then  we first claim that for any $x\in B_2(0)$,
\begin{equation}\label{re:324}
  \varepsilon_{\lambda}^{-p}g(\varepsilon_{\lambda}x+x_{\lambda})\to \infty\,\ \hbox{as\ $\lambda\to\infty$}.
\end{equation}
By the assumption $\rm(M_1)$, it is clear that the claim (\ref{re:324}) is true for the case $x_{\lambda}\nrightarrow 0$ as $\lam\to\infty$.
Therefore, the rest is to consider the case $x_{\lambda}\to 0$ as $\lam\to\infty$. In fact, in the latter case we obtain that $|\varepsilon_{\lambda}x+x_{\lambda}|$ is small in $B_2(0)$. By the assumption $\rm(M_2)$ for $p>0$, we derive from (\ref{re:32T}) that
\begin{align*}
\varepsilon_{\lambda}^{-p}g(\varepsilon_{\lambda}x+x_{\lambda}) &\geq \frac{1}{2}\varepsilon_{\lambda}^{-p}|\varepsilon_{\lambda}x+x_{\lambda}|^p\geq \Big(\frac{1}{2}\Big)^{p+1}\varepsilon_{\lambda}^{-p}|x_{\lambda}|^p\to\infty\,\ \hbox{as}\,\ \lambda\to\infty,
\end{align*}
where $x\in B_2(0)$. This establishes the claim \eqref{re:324}.

We now observe from \eqref{prop u1} that $\varepsilon_{\lambda}^2\mu_{\lambda}\leq C(p)$ as $\lam\to\infty$. Together with \eqref{re:324}, this estimate yields that $\varepsilon_{\lambda}^{-p}g(\varepsilon_{\lambda}x+x_{\lambda})-\varepsilon_{\lambda}^2\mu_{\lambda}>0$ holds for large $\lambda>0$.
Thus, we note from (\ref{eq:Euler}) that for large $\lam >0$,
\begin{equation}\label{eq:324}
-\bigtriangleup w_{\lambda}(x)-c_\lambda(x)w_{\lambda}(x)\leq0\ \ \hbox{in\ $B_2(0)$,\ where \ $c_\lambda(x)=\int_{\R^d}\frac{w^2_\lambda(y)}{|x-y|^2}dy$}.
\end{equation}
By the Hardy-Littlewood-Sobolev inequality (c.f. \cite[Theorem 4.3]{EhL1}), we have
\begin{equation}\label{eq:340}
\|c_\lambda(x)\|_{L^d(B_2(0))}\leq\|c_\lambda(x)\|_{L^d(\R^d)}\leq C\|w_\lambda(x)\|^2_{L^\frac{2d}{d-1}(\R^d)}\leq C_1(N,p)\ \ \hbox{as}\,\ \lambda\to\infty.
\end{equation}
Using the De Giorgi-Nash-Moser theory (c.f. \cite[Theorem 4.1]{HL}), we then obtain that
\begin{equation}\label{eq:346}
w_\lambda(0)\leq C_1(N,p)\Big(\int_{B_{2}(0)}w^2_\lambda(y)dy\Big)^{\frac{1}{2}}\leq C_2(N,p)\ \ \hbox{as}\,\ \lambda\to\infty.
\end{equation}
On the other hand, by Hardy's inequality, we reduce from (\ref{def:317})
\begin{equation}\label{eq:347}
\begin{split}
c_\lambda(0)&=\int_{\R^d}\frac{w^2_\lambda(y)}{|y|^2}dy\leq C(N) \int_{\R^d}|\nabla w_\lambda|^{2}dy\leq C(N,p).
\end{split}
\end{equation}
Since $w_\lambda(x)$ attains its maximum at $x=0$, we obtain from (\ref{eq:Euler}) and (\ref{eq:347}) that
\begin{equation}\label{re:337}
  \varepsilon_{\lambda}^{-p}g(x_{\lambda})w_\lambda(0) \leq \varepsilon_{\lambda}^2\mu_{\lambda} w_\lambda(0)+c_\lambda(0)w_\lambda(0)   \leq\varepsilon_{\lambda}^2\mu_{\lambda} w_\lambda(0)+C(N,p)w_\lambda(0).
\end{equation}
Recall from \eqref{prop u2} that there exists some constant $C_0(N,p)>0$ such that
\begin{equation}\label{eq:337}
w_\lambda(0)\geq C_0(N,p)>0 \ \ \hbox{as}\,\ \lambda\rightarrow\infty.
\end{equation}
By \eqref{eq:346}, (\ref{re:337}) and (\ref{eq:337}), we thus conclude from (\ref{prop u1}) that
\begin{equation}\label{eq:338}
C(p) +C(N,p)
  \geq \varepsilon_{\lambda}^{-p}g(x_{\lambda}) \ \ \hbox{as}\,\ \lambda\rightarrow\infty.
\end{equation}
Since it follows from \eqref{re:324} that  $\varepsilon_{\lambda}^{-p}g(x_{\lambda})\to\infty $ as $\lambda\rightarrow\infty$, we reduce a contradiction from the above estimate (\ref{eq:338}). Therefore, the estimate \eqref{point extimae} holds true.
\vskip 0.1truein

\noindent{(2).} We note from \eqref{eq:Euler} that $w_\lambda(x)$ satisfies
\begin{equation*}
  -\triangle w_\lambda+\big[\varepsilon_{\lambda}^{-p}g(\varepsilon_{\lambda}x+x_{\lambda})-\varepsilon^2_{\lambda}\mu_{\lambda}-
  c_{\lambda}(x)\big]w_{\lambda}=0\,\  \hbox{in}\,\ \R^d,
\end{equation*}
where $c_\lambda(x)=\int_{\R^d}\frac{w^2_\lambda(y)}{|x-y|^2}dy$. If the estimates
\begin{equation}\label{re:3351}
  w_{\lambda}(x)\leq C(N,p)\quad\hbox{in $B_R^c(0)$,}
\end{equation}
and
\begin{equation}\label{re:335}
\varepsilon_{\lambda}^{-p}g(\varepsilon_{\lambda}x+x_{\lambda})-\varepsilon^2_{\lambda}\mu_{\lambda}-c_{\lambda}(x)\geq 1\ \  \hbox{in  $B_R^c(0)$,}
\end{equation}
hold true for sufficiently large $R>0$ and $\lam >0$, it then yields from the comparison principle that \eqref{re exp} holds, see \cite{GWZZ} for the related argument. So, the rest is to establish (\ref{re:3351}) and (\ref{re:335}).

By Lemma \ref{prop u} (1), we have $\varepsilon^2_{\lambda}\mu_{\lambda}\leq C(p)$ as $\lambda\rightarrow\infty$.
We claim that for any $k>0$, there exists a sufficiently large constant $R=R(k)>0$ such that
\begin{equation}\label{re:3352}
 \varepsilon^{-p}_{\lambda}g(\varepsilon_{\lambda}x+x_{\lambda})>k\quad\hbox{in\ $B_R^c(0)$},
\end{equation}
if $\lam >0$ is large enough. In order to prove the claim (\ref{re:3352}), consider a small fixed constant $\varepsilon_0>0$. For any $k>0$, if $|x|>\frac{\varepsilon_0}{\varepsilon_{\lambda}}$, we then have
$$\varepsilon^{-p}_{\lambda}g(\varepsilon_{\lambda}x+x_{\lambda})\geq \varepsilon_\lam ^{-p}g\big(\frac{\varepsilon_{\lambda}x}{2}\big)>k\ \ \hbox{as}\,\  \lambda\rightarrow\infty,
$$
and hence the estimate (\ref{re:3352}) is proved if $|x|>\frac{\varepsilon_0}{\varepsilon_{\lambda}}$.
Next, if $R<|x|\leq \frac{\varepsilon_0}{\varepsilon_{\lambda}}$, we reduce that
$|\varepsilon_{\lambda}x+x_{\lambda}|\leq 2\varepsilon_0$. This further implies from the assumption $(M_2)$ that
$g(\varepsilon_{\lambda}x+x_{\lambda})\geq \frac{1}{2}|\varepsilon_{\lambda}x+x_{\lambda}|^p$ for $R<|x|\leq \frac{\varepsilon_0}{\varepsilon_{\lambda}}$. Hence, for any $k>0$, take $R=\max\{2^{\frac{p+1}{p}}k^{\frac{1}{p}}, 2C(N, p)\}>0$, where $C(N, p)>0$ is as in (\ref{point extimae}), so that for large $\lam>0$,
$$|\varepsilon_{\lambda}x+x_{\lambda}|\geq \frac{1}{2}|\varepsilon_{\lambda}x| \ \ \mbox{for} \,\ R<|x|\leq \frac{\varepsilon_0}{\varepsilon_{\lambda}},$$
which then gives that for large $\lam>0$,
\begin{align*}
\varepsilon_{\lambda}^{-p}g(\varepsilon_{\lambda}x+x_{\lambda})&
\geq \frac{1}{2}\varepsilon_{\lambda}^{-p}|\varepsilon_{\lambda}x+x_{\lambda}|^p
\geq \Big(\frac{1}{2}\Big)^{p+1}\varepsilon_{\lambda}^{-p}| \varepsilon_{\lambda}x|^p\\
&> \Big(\frac{1}{2}\Big)^{p+1}R^p>k \ \ \mbox{for} \,\ R<|x|\leq \frac{\varepsilon_0}{\varepsilon_{\lambda}}.
\end{align*}
Therefore, the claim \eqref{re:3352} is established.

Applying (\ref{prop u1}) and \eqref{re:3352}, we note that there  exists a sufficiently large constant $R=R(N, p)>0$ such that for large $\lam >0$, $$\varepsilon^{-p}_{\lambda}g(\varepsilon_{\lambda}x+x_{\lambda})-\varepsilon^{2}\mu_{\lambda}\geq 0\ \ \hbox{in}\,\ B^c_R(0),$$
and hence
\begin{equation}\label{re:3352E}
-\bigtriangleup w_{\lambda}(x)-c_\lambda(x)w_{\lambda}(x)\le 0\ \ \hbox{in\ $B^c_R(0)$,\ where $c_\lambda(x)=\int_{\R^d}\frac{w^2_\lambda(y)}{|x-y|^2}dy$},
\end{equation}
if $R>0$ and $\lam>0$ are large enough. Similar to proving (\ref{eq:346}), the estimate \eqref{re:3351} then follows by applying the De Giorgi-Nash-Moser theory to (\ref{re:3352E}). Therefore, the estimate \eqref{re:3351} is now proved.

As for \eqref{re:335}, we derive from \eqref{re:3351} that for large $R>0$ and $\lam >0$,
\begin{equation}\label{re:3350}
  \begin{split}
c_{\lambda}(x)&=\int_{\R^d}\frac{w^2_\lambda(y)}{|x-y|^2}dy=\int_{B_1(x)}\frac{w^2_\lambda(y)}{|x-y|^2}dy+\int_{B_1^{c}(x)}\frac{w^2_\lambda(y)}{|x-y|^2}dy\\
              &\leq C(N,p) \int_{B_1(x)}\frac{1}{|x-y|^2}dy+\int_{B_1^{c}(x)}{w^2_\lambda(y)}dy\\
              &\leq C(N,p)\ \ \hbox{in}\,\  B^c_R(0) .
  \end{split}
\end{equation}
By taking sufficiently large $R=R(N, p)>0$, the estimate \eqref{re:335} therefore follows in view of (\ref{prop u1}), \eqref{re:3352} and (\ref{re:3350}).

\vskip 0.1truein

\noindent{(3).}
Following \eqref{prop u2} and (\ref{def:317}), there exist a subsequence  $\{w_k\}$ of $\{w_{\lambda_k}\}$, where $\lambda_k\to\infty$ as $k\to\infty$, and $0\le w_0\in H^1(\R^d)$ such that
\begin{equation}\label{def:318}
 w_k\rightharpoonup w_0\neq 0 \ \ \hbox{in}\,\ H^1(\R^d)\ \ \hbox{as}\,\  k \rightarrow \infty,
\end{equation}
and hence the estimate (\ref{re conver}) holds true.

We next prove (\ref{re conver*}) as follows.
Indeed, for any $l>0$, we have
\begin{equation}\label{eq:341}
\int_{\{g(x)\geq l\varepsilon_{\lambda}^p\}}u^2_\lambda(x)dx\leq
\frac{1}{l\varepsilon_{\lambda}^p}\int_{\R^d}g(x)u^2_\lambda(x)dx\leq\frac{M_2}{l}\ \ \hbox{as}\,\  \lambda\rightarrow\infty.
\end{equation}
We thus deduce from above that for any large $l>0,$
\begin{equation}\label{eq:342}
\begin{split}
\frac{M_2}{l}&\geq\int_{\{g(x)\geq l\varepsilon_{\lambda}^p\}}u^2_\lambda(x)dx=N-\int_{\{g(x)\leq l\varepsilon_{\lambda}^p\}}u^2_\lambda(x)dx \\
             &\geq N-\int_{\{|x|\leq 2l^{\frac{1}{p}}\varepsilon_{\lambda}\}}u^2_\lambda(x)dx=N-\int_{\{|x-\frac{x_\lambda}{\varepsilon_{\lambda}}|\leq 2l^{\frac{1}{p}}\}}w^2_\lambda(x)dx \\
             &\geq N-\int_{\{|x|\leq 4l^{\frac{1}{p}}\}}w^2_\lambda(x)dx\ \ \hbox{as}\,\  \lambda\rightarrow\infty,
\end{split}
\end{equation}
where (\ref{point extimae}) is used in the last inequality.
By (\ref{def:318}), we obtain that for any large $l>0,$
\begin{equation}\label{eq:343}
w_k\rightarrow w_0\,\  \hbox{strongly in}\,\ L^q\big(\{|x|\leq 4l^{\frac{1}{p}}\}\big) \ \ \hbox{as}\,\  k\rightarrow\infty,
\end{equation}
where $2\leq q<2^*$. Therefore, we conclude from (\ref{eq:342}) and (\ref{eq:343}) that for any large $l>0,$
\begin{equation}\label{eq:344}
N\geq\int_{\{|x|\leq 4l^{\frac{1}{p}}\}}w^2_0(x)dx\geq N-\frac{M_2}{l}.
\end{equation}
Since $l>0$ is arbitrary, we conclude from (\ref{eq:344}) that
\begin{equation}\label{eq:345}
N=\|w_0\|^2_2\,\ \hbox{and}\,\ w_k\rightarrow w_0\,\  \hbox{strongly in}\,\ L^2(\R^d)\ \ \hbox{as}\,\  k\rightarrow\infty.
\end{equation}
Because $w_k$ is bounded uniformly in $H^1(\R^d)$, the estimate (\ref{re conver*}) is now proved.
This completes the proof of Lemma \ref{exp decaey}.
\qed

\vskip 0.1truein

\noindent{\bf Proof of Theorem \ref{thm1.2}.} Suppose $u_{\lambda}$ is a positive minimizer of (\ref{def:Ea}), and let  $w_\lambda$ be defined by (\ref{def:316}). Note from Lemma \ref{exp decaey} that $w_\lambda$ satisfies the exponential decay (\ref{re exp}) and there exist a subsequence  $\{w_k\}$ of $\{w_{\lambda_k}\}$, where $\lambda_k\to\infty$ as $k\to\infty$, and $0\le w_0\in H^1(\R^d)$ such that $w_k \rightharpoonup w_0 $ in $ H^1(\R^d) $ and $
w_k\to w_0$ strongly in $L^q(\R^d)$ as $k\to\infty$, where $2\leq q<2^*$. It yields from (\ref{re exp}) that $w_0(x)\leq C(N, p)e^{-|x|}$ in $B_R^c(0)$ for sufficiently large $R>0$. Applying (\ref{point extimae}), up to a subsequence of $\{\lam _k\}$, we have $\lim_{k\to\infty}\frac{x_{\lambda _k}}{\varepsilon_{k}}=y_0$ for some $y_0\in\R^d$, where $x_{\lambda _k}$ is a maximal point of $u_{\lambda _k}$ and $\varepsilon_k=\lam _k^{-\frac{1}{2+p}}>0$. Direct calculations now yield from (\ref{def:316}) that
\[ \begin{split}
  E_{\lambda _k}(u_{\lambda _k}) =& \inte|\nabla u_{\lambda _k}(x)|^2dx+{\lambda _k} \inte g(x)u_{\lambda _k}^2(x)dx-\frac{1}{2} \inte\inte \frac{u_{\lambda _k}^2(x)u_{\lambda_k}^2(y)}{|x-y|^2}dxdy\\
  =&\varepsilon_k^{-2}\inte|\nabla w_k(x)|^2dx+\varepsilon_k^{-(2+p)}\inte g(\varepsilon_kx+x_{\lambda_k})w_k^2(x)dx\\
&-\frac{1}{2}\varepsilon_k^{-2} \inte\inte \frac{w_k^2(x)w_k^2(y)}{|x-y|^2}dxdy\\
  =&  \varepsilon_k^{-2}\inte|\nabla w_k(x)|^2dx+\varepsilon_k^{-2}\inte \Big|x+\frac{x_{\lambda_k}}{\varepsilon_k}\Big|^pw_k^2(x)dx\\
  &-\frac{1}{2}\varepsilon_k^{-2} \inte\inte \frac{w_k^2(x)w_k^2(y)}{|x-y|^2}dxdy\\
   &+\varepsilon_k^{-2}\inte \Big(\varepsilon_k^{-p}g(\varepsilon_kx+x_{\lambda_k})-
   \Big|x+\frac{x_k}{\varepsilon_k}\Big|^p\Big)
   w_k^2(x)dx.
 \end{split}
\]
Similar to (\ref{eq:gx}), one can prove that
\[
\inte \Big(\varepsilon_k^{-p}g(\varepsilon_kx+x_{\lambda_k})-
   \Big|x+\frac{x_k}{\varepsilon_k}\Big|^p\Big)   w_k^2(x)dx=o(1)\ \ \hbox{as}\,\  k\rightarrow\infty.
\]
We then deduce from above that
\begin{equation}\label{re last}
 \begin{split}
&\quad\liminf_{k \to \infty}{\lambda _k}^{-\frac{2}{2+p}}e_{\lambda_k}(N)=\liminf_{k \to \infty}\varepsilon_k^{2}E_{\lambda_k}(u_{\lambda _k})\\
\geq &\inte|\nabla w_{0}(x)|^2dx+\inte |x+y_0|^pw_{0}^2(x)dx-\frac{1}{2} \inte\inte \frac{w_{0}^2(x)w_{0}^2(y)}{|x-y|^2}dxdy\\
=&\inte|\nabla \hat {w}_{0}(x)|^2dx+\inte |x|^p\hat{w}_{0}^2(x)dx-\frac{1}{2} \inte\inte \frac{\hat{w}_{0}^2(x)\hat{w}_{0}^2(y)}{|x-y|^2}dxdy
\geq e_{\infty}(N),
 \end{split}
\end{equation}
 where $\hat w_{0}(x)=w_{0}(x-y_0)\ge 0$. Following Lemma \ref{lem:energy estimate}, the above inequalities must be identities,which further imply that $\hat{w}_0\ge 0$ is a  minimizer of $e_{\infty}(N)$. Therefore,
$$e_{\infty}(N)=E_{\infty}(\hat{w}_0(x)),\ \ 0<N<N^*,$$
and
\begin{equation}\label{re lastA}
\lim_{k\to\infty}\|\nabla w_k(x)\|_2^2=\|\nabla  w_0(x)\|_2^2\ \ \mbox{and}\,\ w_k(x)\rightarrow w_0(x)=\hat w_{0}(x+y_0) \,\ \hbox{in}\,\ H^1(\R^d)\ \ \mbox{as}\,\ k\to\infty.
\end{equation}
which implies that the limit behavior (\ref{def:18}) now follows from (\ref{def:316}) and (\ref{re lastA}).
Also, $\hat{w}_0\ge 0$ satisfies the Euler-Lagrange equation (\ref{def:333}) for some Lagrange multiplier $\mu =\mu(N)\in\R$ such that $\int_{\R^d} \hat w^2_0=N$. Moreover, the maximum principle applied to (\ref{def:333}) yields that $\hat{w}_0>0$ in $\R^d$. Similar to \cite{GWZZ}, one can further derive that $\hat{w}_0>0$ is essentially a ground state of (\ref{def:333}).

The rest is to prove that $y_0=0$. On the contrary, suppose that $|y_0|>0$. Since $w_k(0)= \max_{x\in\R^d}w_k(x)$ for all $k>0$, we deduce from (\ref{re lastA}) that $\hat{w}_0(|y_0|)=\max_{x\in\R^d}\hat w_0(x)$, which is a contradiction to the fact that $\hat{w}_0>0$ is strictly decreasing in $|x|$ from Proposition \ref{prop:3.1}.
Therefore, we have $y_0=0$, which then completes the proof of Theorem 1.2.
\qed

\vspace {.5cm}

\noindent {\bf Acknowledgements:} The first author thanks Prof. Robert Seiringer for introducing him the related interesting phenomena addressed in this paper and for constant encouragements over the past few years. Part of this paper was finished when the first author visited Pacific Institute for Mathematical Sciences (PIMS) at UBC from January to February in
2018. He would like to thank PIMS for their warm hospitality.

\end{document}